\documentclass[11pt]{article}
\usepackage[utf8]{inputenc}
\usepackage[affil-it]{authblk} %
\usepackage{cite}
\usepackage{setspace}
\usepackage[title]{appendix}
\usepackage{mathrsfs}
\usepackage{graphicx}
\usepackage{multirow}
\usepackage[labelsep=none]{caption}
\usepackage {enumerate}
\usepackage{float}   
\usepackage{subfigure} 

\usepackage{amsfonts}
\usepackage{amssymb}
\usepackage{amsmath}
\usepackage{amsthm}
\usepackage[colorlinks=true]{hyperref}
\hypersetup{linkcolor=blue}
\theoremstyle{plain}
\newtheorem{thm}{Theorem}[section]

\newtheorem{lem}[thm]{Lemma}

\newtheorem*{remark}{Remark}  
\theoremstyle{definition}

\newcommand{\intend}[1]{\noindent #1}
\numberwithin{equation}{section}

\date{}

\begin{document}
	\begin{spacing}{0.5}
	\end{spacing}
	\title{{\Large\bf {Asymptotic Behavior of Rupture Solutions for the Elliptic MEMS Equation with Hénon and External Pressure Terms}}}
	\author[1]{Yunxiao Li}
    \author[2]{Yanyan Zhang}

	\affil[1]{{\footnotesize School of Mathematical Sciences, East China Normal University, Shanghai 200241, People's Republic of China}}
   \affil[2]{{\footnotesize School of Mathematical Sciences, Key Laboratory of MEA (Ministry of Education) and Shanghai Key Laboratory of PMMP, East China Normal University, Shanghai 200241, China}}

		\maketitle
	\begin{minipage}{14cm} {\bf Abstract:} This paper investigates an elliptic MEMS-Type equation with Hénon and external pressure terms:
 \begin{equation*}
\begin{cases} 
\Delta u = \dfrac{\lambda |x|^{\alpha}}{u^{p}}+F  & x \in \mathbb{R}^{N} \setminus \{0\}, \\[4pt]
u(0) = 0 ~\text{and} ~u > 0  & x \in \mathbb{R}^{N} \setminus \{0\} , \\
\end{cases}
 \end{equation*}
where $ N \ge 1, \lambda > 0, p > 0, \alpha > -2$ and $F \in R $ are constants. We study positive rupture solutions with rupture point at the origin ($u(0)=0$). Our main emphasis is on asymptotic radial rupture solutions: we prove the existence of both radial and non-radial solutions, characterize their asymptotic behavior near the origin, and obtain a full asymptotic expansion of arbitrary order

  
	\end{minipage}

\section{Introduction}
In this paper, we investigate the following elliptic MEMS-type equation with Hénon term:
\begin{equation}
\tag{1}\label{QQ}
\begin{cases} 
\Delta u = \dfrac{\lambda |x|^{\alpha}}{u^{p}}+F  & x \in \mathbb{R}^{N} \setminus \{0\}, \\[4pt]
u(0) = 0 ~\text{and}~ u > 0  & x \in \mathbb{R}^{N} \setminus \{0\} , \\
\end{cases}
\end{equation}
where $ N \ge 1, \lambda > 0, p > 0, \alpha > -2$ and $F \in R $ are constants. We establish the existence of solutions to \eqref{QQ} and derive asymptotic expansions in a sufficiently small neighborhood of the origin.

Micro-Electro-Mechanical Systems (MEMS) are devices that integrate miniature mechanical components with electronics using microfabrication technology. We consider a canonical MEMS model, as illustrated in Figure 1 (c.f.\cite{24z}), consisting of a deformable elastic membrane suspended above a fixed ground plate. When a voltage is applied, the electrostatic force attracts the membrane towards the plate. As the voltage increases, the membrane may eventually touch the plate. This phenomenon, known as "pull-in instability", can compromise the stable operation of high-precision microelectronic devices. Conversely, in certain applications such as airbag restraint systems or inkjet printers, such contact is desirable.

In this work, we focus on the profile of the membrane near the touchdown point—referred to mathematically as the rupture point—at the moment contact occurs. A systematic characterization of the local asymptotic profile near the rupture point, along with its dependence on system parameters, provides essential guidance for the design and reliability assessment of MEMS devices.

\begin{figure}[htb]
	\centering
	\includegraphics[width=0.6\textwidth,trim=3.5cm 0cm 0cm 0cm,clip]{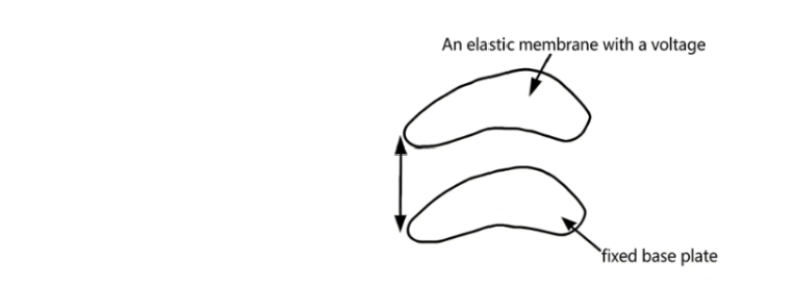}
  \caption{}
\end{figure}

Mathematically, let $u \ge 0$ denote the distance between the membrane and the ground plate. Let $\lambda > 0$ represent the applied voltage, $|x|^\alpha~(\alpha > -2)$ represent the variable permittivity profile of the membrane, and $F$ denote the external pressure (where $F > 0$ corresponds to a downward force and $F < 0$ to an upward force). We focus on the case where the membrane ruptures at $x=0$. This is modeled by the equation (cf. \cite{23z,24z,25z}):
\begin{equation*}
\begin{cases} 
\Delta u = \dfrac{\lambda |x|^{\alpha}}{u^{p}}+F  & x \in \mathbb{R}^{2} \setminus \{0\}, \\[4pt]
u(0) = 0 ~\text{and}~ u > 0  & x \in \mathbb{R}^{2} \setminus \{0\} . \\
\end{cases}
\end{equation*}
Problem \eqref{QQ} is obtained by generalizing the exponent $p=2$ to a general $p>0$, extending the domain to the N-dimensional setting ($N \ge 1$).

Problem \eqref{QQ}  has been studied in various parameter regimes. For $\alpha=0$, the rupture solutions of problem \eqref{QQ} were investigated in \cite{19z,21z}. In particular, for $\lambda=1$, $p>1$ and $N\geq2$, the Hölder continuity and the Hausdorff dimensions of rupture sets were studied in \cite{21z}. For $\lambda=1$, $F=0$ and $N\geq3$, \cite{19z} obtained infinitely many non-radial rupture solutions for $p > \frac{N-1-2\sqrt{N-2}}{2\sqrt{N-2}-(N-5)}$.  For $\alpha \neq 0$, the rupture solutions of problem \eqref{QQ} were investigated in \cite{7z,10z,12z}. In particular, for $N=2$, $p=2$ and $F \ge 0$, \cite{7z} investigated the isotropic and anisotropic singularities of rupture solutions with respect to different values of $\alpha$. For $N=2$ and $F=0$, \cite{10z} analyzed the isotropic and anisotropic behaviors in various parameter regions of $(\alpha, p)$ and studied global solutions on $\mathbb{R}^2 \setminus \{0\}$. For $\lambda=1$ and $F=0$, a dichotomy result regarding the Morse index of such solutions was presented in \cite{12z}, depending on whether the exponent $p$ is below or above the critical threshold .

Significant progress has also been made regarding the asymptotic expansion of solutions near singular points. For instance, \cite{7z} derived the first two terms of the expansion for asymptotically isotropic rupture solutions. Using spherical harmonics, \cite{4z}  derived the asymptotic expansion near an isolated singularity for the Yamabe equation $-\Delta u = \frac{1}{4}n(n-2)u^{\frac{n+2}{n-2}}$. More recently, \cite{8z} utilized spherical harmonics to analyze asymptotic expansions to arbitrary order for the steady thin-film-type equation near the origin. This corresponds to a special case of Equation \eqref{QQ} with $\lambda=1$ and $\alpha=0$.

In this paper, we generalize these analyses to the general $\alpha$. Regarding asymptotic radial solutions to problem \eqref{QQ}, we investigate the existence and asymptotic behavior of radial and non-radial rupture solutions near the origin. An asymptotic radial solution is defined as a solution of the form$$u(x)=C|x|^d\left(1+O\left(|x|^\epsilon\right)\right) \quad as ~|x| \rightarrow 0^{+}$$
for some $\epsilon>0$, where $C$ and $d$ are non-zero constants. The main results are stated as follows.

\begin{thm}  

For the elliptic MEMS problem \eqref{QQ} with $N\geq2$, $F\neq0$ and $-2 <\alpha < 2p $, there exists at least one radial rupture solution near the origin satisfying

$$
u(r)=\Lambda r^{\frac{\alpha+2}{p+1}}+\sum_{i=1}^\infty d_i r^{\frac{(2 p-\alpha) i+(\alpha+2)}{p+1}}=\Lambda r^{\frac{\alpha+2}{p+1}}\left(1+O\left(r^{\frac{2 p-\alpha}{p+1}}\right)\right) \qquad \text{as } r \rightarrow 0^{+},
$$
where $d_1$ depends on $F,p,\alpha,N$ and $d_2,\dots,d_\infty$ depend on $p,\alpha,N,\lambda$.
\begin{equation}
\tag{2}\label{qcc}
\Lambda=\left(\frac{\lambda}{\frac{\alpha+2}{p+1}\left(\frac{\alpha+2}{p+1}+N-2\right)}\right) ^{\frac{1}{p+1}}.
\end{equation}
\end{thm}

\begin{remark}
 For the case of $F=0$, it has been proved in \cite{12z} that $u(x)=\Lambda |x|^{\frac{\alpha+2}{p+1}}$ is the only positive radial rupture solution of \eqref{QQ}.
\end{remark}

\begin{thm}
For the elliptic MEMS problem \eqref{QQ} with $N\geq2$, $F=0$ and $\alpha>-2  $, or $N\geq2$, $F\neq0$ and $-2<\alpha<2p  $, there exist infinitely many non-radial positive asymptotic radial rupture solutions satisfying
$$
u(x)=\Lambda |x|^{\frac{\alpha+2}{p+1}}+\sum_{j=1}^\infty \sum_{i=0}^{j-1} C_{j i}\left(\frac{x}{|x|}\right)(\ln |x|)^i|x|^{{\mu_j}+{\frac{\alpha+2}{p+1}}}=\Lambda |x|^{\frac{\alpha+2}{p+1}}\left(1+O\left(|x|^{\mu_1}\right)\right)
$$
as $|x| \rightarrow 0^{+}$,
where $\Lambda$ is from \eqref{qcc}, $\left\{\mu_j\right\}$ is a strictly increasing sequence of positive numbers diverging to $+\infty$, and
$$
{\mu_1}=
\sigma_1^{(k)}  ~~~ (\text{if}~~\delta^{(k-1)}<\alpha<\delta^{(k)}, k=1,2,\dots)
$$
\intend{in the case $F=0$.}

$$\mu_1=\begin{cases}
{\sigma_1^{(k)}}& \text{if~ $\delta^{(k-1)}<\alpha<\delta^{(k)}$ and $\alpha<2p-(p+1)\sigma_1^{(k)}$, } k=1,2,\dots\\
\displaystyle\frac{2 p-\alpha}{p+1} &\text{if ~$\delta^{(k-1)}<\alpha<\delta^{(k)}$ and $\alpha\geq2p-(p+1)\sigma_1^{(k)}$, }k=1,2,\dots
\end{cases}$$
\intend{in the case $F\neq0$, where}
$$
\delta^{(k)}=\begin{cases}
-2 & \text{if $k=0,$}\\
\frac{-1}{2} (p+1)(N+2)+\frac{1}{2}\sqrt{(N-2)^2(p+1)^2+4(p+1) k(N-2+k)} & \text{if $k\geq 1 $}
\end{cases}
$$
and

$\sigma_1^{(k)} =-\frac{1}{2}\left( N - 2 + 2\frac{\alpha+2}{p+1} \right) + \frac{1}{2} \sqrt{ \left( N - 2 + 2\frac{\alpha+2}{p+1}  \right)^2 + 4k(N - 2 + k) - 4(\alpha+2)\left( N - 2 + \frac{\alpha+2}{p+1}  \right) }.$
\end{thm}
\begin{remark}
 For the case of $N=1$, there exists at least one solution near the origin satisfying
\[
u(x)=\Lambda^{'} |x|^{\frac{\alpha+2}{p+1}}+\sum_{i=1}^\infty d_i x^{\frac{(2 p-\alpha) i+(\alpha+2)}{p+1}}=\Lambda^{'} |x|^{\frac{\alpha+2}{p+1}}\left(1+O\left(x^{\frac{2 p-\alpha}{p+1}}\right)\right) \qquad \text{as } x \rightarrow 0^{+},
\]
where $\Lambda^{'}=\left(\frac{\lambda}{\frac{\alpha+2}{p+1}\left(\frac{\alpha+2}{p+1}-1\right)}\right) ^{\frac{1}{p+1}}$, since the arguments for solution in the case $N=1$ is completely similar to Theorem 1.1, we shall omit the case $N=1$.
\end{remark}
\vspace{2em}

To facilitate the proof of the subsequent theorems, we introduce the following change of variables. Let $t=\ln |x|$, $\theta=\dfrac{x}{|x|}$, then
\begin{equation}
\tag{3}\label{Qz}
\Delta u=u_{r r}+\frac{N-1}{r} u_r+\frac{1}{r^2} \Delta_\theta u=\frac{\lambda|x|^\alpha}{u^p}+F.
\end{equation}

Let $z(t, \theta)=r^{-\frac{\alpha+2}{p+1} }u(x)-\Lambda$, where $\Lambda$ is from \eqref{qcc}, then we obtain:
\begin{equation}
\tag{4}\label{3}
z_{t t}+\left(N-2+2 {\frac{\alpha+2}{p+1}}\right) z_t+(\alpha+2)(N-2+{\frac{\alpha+2}{p+1}}) z+\Delta{ }_\theta z=f(z) +F  e^{\frac{2 p-\alpha}{p+1} t}  ,
\end{equation}
where $f(z) = \dfrac{\lambda}{(z+\Lambda)^p}-\dfrac{\lambda}{\Lambda^p}+\dfrac{\lambda p z}{\Lambda^{p+1}}$.

\vspace{1em}

We highlight the key innovations and distinguishing features of our study compared to the existing literature as follows:

\vspace{1em}

\textbf{(A) Extension to General Parameter Regimes.} Existing literature has primarily addressed equation \eqref{QQ} within limited parameter regimes, such as the specific case $N=2, p=2$ studied in [8], or the case assuming $\lambda=1$ and $\alpha=0$ in [9]. We generalize these analyses to the framework of general $\alpha$. This extension is nontrivial as it introduces significant analytical challenges not present in the isotropic or restricted parameter cases.

\vspace{1em}

\textbf{(B) Overcoming Spectral Difficulties Involving Complex Eigenvalues.} A major technical obstacle in the general $\alpha$ setting is that the associated linear operator may possess complex eigenvalues, whereas previous works mostly dealt with real spectrums. We overcome this difficulty by performing a more involved classification and establishing delicate estimates.

\vspace{1em}

\textbf{(C) Theoretical foundation for the design of MEMS models.} While asymptotic expansions have been studied for simpler cases, our work establishes these precise characterizations within a general framework relevant to broader physical contexts. These results provide essential theoretical insights for the optimization of MEMS device structures, allowing for a more accurate description of rupture profiles.

\vspace{1em}

The remainder of this paper is organized as follows. In Section 2, we derive asymptotic expansions of arbitrary order near the origin for the radial case. Section 3 extends these results non-radial case. Finally, in Section 4, we present the proofs of Theorem 1.1 and 1.2.

\section{Asymptotic behavior for radial solution}
In this section, we give the arbitrary order asymptotic expansion of radial case.
\begin{thm}
    
Suppose that $u(x)$ is a positive radial rupture solution of \eqref{QQ} with $F \neq 0$. Assume it exists $\varepsilon>0$ such that
\[
u(r)=\Lambda r^{\frac{\alpha+2}{p+1}}\left(1+O\left(r^{\varepsilon}\right)\right)~~
as~ r \rightarrow 0^{+},
\]
where $\Lambda$ is the same as in theorem 1.1, set
$$
  t=\ln r, z(t)=r^{-\frac{\alpha+2}{p+1}} u(r)-\Lambda, 
$$
then for any $k \gg 1$ we have
\begin{equation*}
z(t) = \sum_{\ell=1}^{k} c_{\ell} e^{\ell \rho t} + O \left( e^{(k+1)\rho t} \right),
\end{equation*}
where $c_\ell$ are constants and
$\rho=\frac{2p-\alpha}{p+1}$.
\end{thm}

\medskip
\begin{proof}
It follows from \eqref{3} that
\[
\left\{
\begin{array}{l}
\displaystyle
z_{t t}+\left(N-2+2\frac{\alpha+2}{p+1}\right) z_t+(\alpha+2)\left(N-2+\frac{\alpha+2}{p+1}\right) z
=f(z)+F e^{\frac{2 p-\alpha}{p+1} t}, \quad t \in(-\infty, 0), \\[6pt]
\displaystyle
z(t)=O\left(e^{\varepsilon t}\right) \rightarrow 0 \quad \text{as } t \rightarrow-\infty,
\end{array}
\right.
\]
where
\[
f(z)=\lambda(\Lambda+z)^{-p}-\lambda \Lambda^{-p}+\lambda p z \Lambda^{-(p+1)}=O(z^2)\quad ~~ as ~z(t) \rightarrow 0.
\]

The associated homogeneous ordinary differential equation is
\[
z_{tt} + \left(N - 2 + 2\frac{\alpha+2}{p+1}\right) z_t
+ (\alpha+2)\left(N-2+\frac{\alpha+2}{p+1}\right) z = 0,
\]
whose characteristic equation is
\[
\sigma^2 + \left(N - 2 + 2\frac{\alpha+2}{p+1}\right) \sigma
+(\alpha+2)\left(N-2+\frac{\alpha+2}{p+1}\right) = 0,
\]
it admits two roots, denoted by $\sigma_1^{(0)}$ and $\sigma_2^{(0)}$.

\noindent\hspace{-0.2em}\textbf{Case1}.\quad $\alpha>\frac{(p+1)(N-2)}{2}(\sqrt{1+\frac{1}{p}}-1)-2,$
\[
\left\{
\begin{aligned}
\sigma_1^{(0)} = -\frac{1}{2} \left(N - 2 + 2\frac{\alpha+2}{p+1} \right) + \frac{i}{2} \sqrt{4(\alpha+2) \left(N - 2 + \frac{\alpha+2}{p+1} \right) - \left(N - 2 + 2\frac{\alpha+2}{p+1} \right)^2 },
\\
\sigma_2^{(0)} = -\frac{1}{2} \left(N - 2 + 2\frac{\alpha+2}{p+1} \right) -  \frac{i}{2} \sqrt{4(\alpha+2) \left(N - 2 + \frac{\alpha+2}{p+1} \right) - \left(N - 2 + 2\frac{\alpha+2}{p+1} \right)^2 },
\end{aligned}
\right.
\]

\vspace{2em}

\noindent\hspace{-0.2em} \textbf{Case2.}\quad $\alpha=\frac{(p+1)(N-2)}{2}(\sqrt{1+\frac{1}{p}}-1)-2,$

\[
\sigma_1^{(0)} = \sigma_2^{(0)} =-\frac{1}{2} \left(N - 2 + 2\frac{\alpha+2}{p+1} \right)  < 0;
\]

\vspace{2em}
\noindent\hspace{-0.2em}\textbf{Case3.}\quad $-2<\alpha<\frac{(p+1)(N-2)}{2}(\sqrt{1+\frac{1}{p}}-1)-2,$

\[
\left\{
\begin{aligned}
\sigma_1^{(0)} = -\frac{1}{2} \left(N - 2 + 2\frac{\alpha+2}{p+1} \right) + \frac{1}{2} \sqrt{ \left(N - 2 + 2\frac{\alpha+2}{p+1} \right)^2-4(\alpha+2) \left(N - 2 + \frac{\alpha+2}{p+1} \right)  }<0,
\\
\sigma_2^{(0)} = -\frac{1}{2} \left(N - 2 + 2\frac{\alpha+2}{p+1} \right) -  \frac{1}{2} \sqrt{ \left(N - 2 + 2\frac{\alpha+2}{p+1} \right)^2 -4(\alpha+2) \left(N - 2 + \frac{\alpha+2}{p+1} \right) }<0,
\end{aligned}
\right.
\]

We now prove
\[
z(t) = O\left(  e^{\frac{2p-\alpha}{p+1}t} \right)\qquad \text{as } t \to -\infty.
\]

We only prove case 3 and other cases are similar. By ordinary equation theory, for $T \ll -1$ and $t\in (-\infty, T)$,
\[
\begin{aligned}
z(t) = & A_1 e^{\sigma_1^{(0)} t}+ A_2 e^{\sigma_2^{(0)} t} 
 + B_1 e^{\sigma_1^{(0)} t} \int_{-\infty}^{t} e^{-\sigma_1^{(0)} s} \left[ f(z(s)) +F e^{\frac{2p-\alpha}{p+1}s} \right]  ds \\ &
+ B_2 e^{\sigma_2^{(0)} t} \int_{-\infty}^{t} e^{-\sigma_2^{(0)}  s} \left[ f(z(s)) +F e^{\frac{2p-\alpha}{p+1}s} \right]  ds,
\end{aligned}
\]

Here $A_1, A_2$ are constants,$B_1, B_2$ only depend on $\sigma_1^{(0)}$ and $\sigma_2^{(0)}$.

We also know
\[
f(z) = O(z^2)=O(e^{2 \epsilon t}).
\]

On the other hand, $z(t) \to 0$ as $t \to -\infty$ ,  so $A_1 = A_2 = 0$,
\begin{equation}
\tag{5}\label{QD}
z(t) = B_1 e^{{\sigma_1^{(0)}}t} \int_{-\infty}^{t} e^{-\sigma_1^{(0)} s} \left[ f(z(s)) +F e^{\frac{2p-\alpha}{p+1}s} \right]  ds  + B_2 e^{{\sigma_2^{(0)}}t} \int_{-\infty}^{t} e^{-\sigma_2^{(0)}  s} \left[ f(z(s)) +F e^{\frac{2p-\alpha}{p+1}s} \right]  ds
\end{equation}

We consider two conditions: $(i)\ 0 < \varepsilon < \frac{2p-\alpha}{p+1}$, $(ii)\ \varepsilon \geq  \frac{2p-\alpha}{p+1}$.

For (ii), $f(z) = O \left( e^{\frac{2p-\alpha}{p+1} t} \right)$,it can be proved by \eqref{QD}.

For (i), $f(z) = O(e^{2 \varepsilon t})$,we know from \eqref{QD} that
\begin{equation}
z(t) = O(e^{2 \varepsilon t}).
\tag{6}\label{GH}
\end{equation}

put \eqref{GH} into $f(z(t)) +F e^{\frac{2p-\alpha}{p+1} t}$ we find
\[
f(z(t)) +F e^{\frac{2p-\alpha}{p+1} t} = O \left( e^{\min \left\{ \frac{2p-\alpha}{p+1}, 4 \varepsilon \right\} t} \right),
\]
we get by \eqref{QD} that:
\begin{equation}
z(t) = O \left( e^{\min \left\{ \frac{2p-\alpha}{p+1}, 4 \varepsilon \right\} t} \right),
\tag{7}\label{6}
\end{equation}

we also consider two conditions $4 \varepsilon \geq \frac{2p-\alpha}{p+1}$ and $4 \varepsilon < \frac{2p-\alpha}{p+1}$ and use arguments similar to the above to obtain conclusions eventually.

Set $\rho={\frac{2p-\alpha}{p+1}}$,

\[
\begin{aligned}
f(z) & =\frac{\lambda \Lambda^{p+1}\left(1-\left(1+\frac{z}{\Lambda}\right)^p+p \frac{z}{\Lambda}\left(1+\frac{z}{\Lambda}\right)^p\right)}{(1+z)^p \Lambda^{p+1}} \\
& =d_2 z^2+d_3 z^3+\cdots+d_n z^n+\cdots \\
& =d_2 e^{2  \rho t}+d_3 e^{3  \rho t}+\cdots+d_n e^{n  \rho t}+\cdots
\end{aligned}
\]

Put in into \eqref{QD} we get that for any $k \gg 1$ and $t \in (-\infty, -1]$
\begin{equation*}
z(t) = \sum_{\ell=1}^{k} c_{\ell} e^{l \rho t} + O \left( e^{(k+1)\rho t} \right).
\end{equation*}
Therefore the proof of Theorem 2.1 is completed.
\end{proof}

\section{Asymptotic behavior for nonradial with asymptotic radial solution }
In this section, we give the arbitrary order asymptotic expansion for the asymptotic radial case.
\begin{thm}
    
Assume $u \in C^2(B \backslash\{0\})$ is a positive rupture solution of \eqref{QQ} and there exists $\epsilon>0$ satisfying $u(x)=\Lambda r^{\frac{\alpha+2}{p+1}}\left(1+O\left(|x|^{\epsilon}\right)\right)$. Let $z(t, \theta)=r^{-\frac{\alpha+2}{p+1}} u(x)-\Lambda$, $\Lambda$ is the same as in Theorem 1.2, $t=\ln |x|$, $\theta=\frac{x}{|x|}$, then there exists a positive sequence $\left\{\mu_j\right\}_{j \geq 1}$ strictly increasing to $\infty$, for any positive integer $n \gg 1$ and $(t, \theta) \in(-\infty,-1) \times S^{N-1},$ 
$$
z(t, \theta)=\sum_{j=1}^n \sum_{l=0}^{j-1} c_{j l}(\theta) t^l e^{\mu_j t}+O\left(t^n e^{\mu_{n+1} t}\right),
$$
where $\sigma_1^{(k)}$ and $\delta^{(k)}$ are the same as in Theorem 1.2,
$c_{jl}(\theta) = \sum_{i=0}^{m_{jl}} a_{jli} Q_i(\theta)$,
$a_{jli}$ is a constant, $m_{jl}$ is a integer depending on $N,j,l,p,a$. $Q_i(\theta)$ is a linear combination of characteristic functions of
\[
- \Delta_{S^{N-1}} Q(\theta) = \lambda_i Q(\theta).
\]
\end{thm}
\vspace{1em}  

It's easy to see that $z(t,\theta) = O(e^{\epsilon t})$ as $t \to -\infty$ and 
\[
z_{tt} + \left( N - 2 + 2\frac{\alpha+2}{p+1} \right) z_t + \Delta_{S^{N-1}} z + (\alpha+2)\left( N - 2 + \frac{\alpha+2}{p+1} \right) z = f(z) +Fe^{\frac{2p-\alpha}{p+1}t} , \quad (t,\theta) \in (-\infty,0) \times S^{N-1},
\]
where
$
 f(z)=\lambda(\Lambda+z)^{-p}-\lambda \Lambda^{-p}+\lambda p z \Lambda^{-(p+1)}=O(z^2) .
$

Considering linearization operator
\begin{equation}
\mathcal{L} = \frac{\partial^2}{\partial t^2} + \left( N - 2 +2\frac{\alpha+2}{p+1}  \right)\frac{\partial}{\partial t} + \Delta_{S^{N-1}} + (\alpha+2)\left( N - 2 + \frac{\alpha+2}{p+1} \right).    
\label{eq:L}\tag{8}
\end{equation}

The operator $\mathcal{L}$ can be divided into infinite  partial operators 
\begin{equation}
\tag{9}\label{wqwqq}
\mathcal{L}_k = \frac{d^2}{dt^2} + \left( N - 2 + +2\frac{\alpha+2}{p+1}\right)\frac{d}{dt} - \lambda_k +  (\alpha+2)\left( N - 2 + \frac{\alpha+2}{p+1} \right),
\end{equation}
for $k = 0,1,2,\ldots$, $\lambda_k$ is the k-th eigenvalues of the eigenvalue problem 
\[
- \Delta_{S^{N-1}} Q = \lambda Q,
\]
 $\lambda_k = k(N-2 + k)$ and
$
m_k = \frac{(N - 2 + 2k)(N - 3 + k)!}{k!(N - 2)!}
$ is the corresponding multiplicity.

We define $\{ Q_1^k(\theta), \ldots, Q_{m_k}^k(\theta) \}$ to be an orthonormal basis corresponding to the eigenspace of  $\lambda_k$ .

The characteristic equation corresponding to \eqref{wqwqq} is
\[
\sigma^2 + \left( N - 2 + 2\frac{\alpha+2}{p+1} \right)\sigma + \left[ (\alpha+2)\left( N - 2 + \frac{\alpha+2}{p+1} \right) - k(N - 2 + k) \right] = 0.
\]

Its two roots are 
{\small
\[
\sigma_1^{(k)} = -\frac{1}{2}\left( N - 2 + 2\frac{\alpha+2}{p+1} \right) + \frac{1}{2} \sqrt{ \left( N - 2 + 2\frac{\alpha+2}{p+1}  \right)^2 + 4k(N - 2 + k) - 4(\alpha+2)\left( N - 2 + \frac{\alpha+2}{p+1}  \right) },
\]
}
{\small
\[
\sigma_2^{(k)} = -\frac{1}{2}\left( N - 2 + 2\frac{\alpha+2}{p+1} \right) - \frac{1}{2} \sqrt{ \left( N - 2 + 2\frac{\alpha+2}{p+1}  \right)^2 + 4k(N - 2 + k) -  4(\alpha+2)\left( N - 2 + \frac{\alpha+2}{p+1}  \right) }.
\]
}

Fix $k=k_0$, so
\[
\sigma_1^{(k_0+1)}>\sigma_1^{(k_0)}>\sigma_1^{(1)}>0, \quad \sigma_2^{(k_0+1)}<\sigma_2^{(k_0)}<\sigma_2^{(1)}<0 
\]
where $-2<\alpha<\delta^{(1)}$ ($k _0\geq2$).

\[
\sigma_1^{(k_0+1)}>\sigma_1^{(k_0)}>0>\sigma_1^{(k_0-1)}
,\quad \sigma_2^{(k_0+1)}<\sigma_2^{(k_0)}<\sigma_2^{(1)}<0 
\]
where $\delta^{(k_0-1)}<\alpha<\delta^{(k_0)}$ ($k_0 \geq2$)

Remark. In fact, it will be possible that $\sigma_1^{(k)}$ and $\sigma_2^{(k)}$ are imaginary numbers when k is sufficiently small. Thew arguments is similar to the proof below because such $\sigma^{(k)}$ is finite and $-\frac{1}{2}\left( N - 2 + 2\frac{\alpha+2}{p+1} \right) <0$.
\vspace{1em}

\begin{lem}
For $N \geq 2$, $p > 0$, $k\geq1$, assume $u$ is a positive rupture solution of \eqref{QQ} satisfying $u(x)=u_s(|x|)\left(1+O\left(|x|^\epsilon\right)\right)$ with
$$
\begin{cases}
0<\epsilon \leq \sigma_1^{(k)}, &\text{as } \delta^{(k-1)}<\alpha<\delta^{(k)} ~~and~~ \alpha<2p-(p+1)\sigma_1^{(k)}, \\
0<\epsilon \leq \frac{2 p-\alpha}{p+1}, &\text{as } \delta^{(k-1)}<\alpha<\delta^{(k)}~~ and ~~ \alpha\geq2p-(p+1)\sigma_1^{(k)} .
\end{cases}
$$
If we define $z(t, \theta) = |x|^{-\frac{\alpha+2}{p+1}} u(x) - \Lambda$, then $z(t, \theta) = \mathcal{O}(e^{\varepsilon t})$ for t $\in (-\infty,1),$ and
$$
\max_{S^{N-1}}|z(t, \theta)| \leq 
\begin{cases}
\sigma_1^{(k)}& \text{as } \delta^{(k-1)}<\alpha<\delta^{(k)} \text{ and } \alpha<2p-(p+1)\sigma_1^{(k)},\\
\frac{2 p-\alpha}{p+1} &\text{as } \delta^{(k-1)}<\alpha<\delta^{(k)} \text{ and } \alpha\geq2p-(p+1)\sigma_1^{(k)} 
\end{cases}
$$
for the case $F\neq0$.
$$
\max_{S^{N-1}}|z(t, \theta)| \leq 
{\sigma_1^{(k)}}
$$
for the case $F=0$.

\end{lem}

\vspace{1em}
\text{Let}
\[
\overline{z}(t) = \frac{\int_{\mathbb{S}^{N-1}} z(t, \theta) \, d\theta}{|\mathbb{S}^{N-1}|}, \quad w(t, \theta) = z(t, \theta) - \overline{z}(t).
\]

Since the Laplace operator on the sphere has zero integral over $\mathbb{S}^{N-1}$, we observe that $w$ satisfies the following equation:                                                                                                                                                                                                                                                                                                                                                                                                                                                                                           
{\small
\begin{equation}
w_{tt} + \left(N - 2 + 2\frac{\alpha+2}{p+1} \right) w_t + \Delta_{\mathbb{S}^{N-1}} w + (\alpha+2) \left( N - 2 + \frac{\alpha+2}{p+1} \right) w = f(z) - \overline{f(z)},  
\tag{10}\label{qcxx}
\end{equation}
}
It is clear that $w(t, \theta) = \mathcal{O}(e^{\varepsilon t})~~as~t \to -\infty$.

A direct computation shows that
\begin{equation}
f(z) - \overline{f(z)} = f'(\xi)\, w - f'(\xi)\, \overline{w},
\tag{11}\label{9}
\end{equation}
where $\xi$ lies between $z$ and $\overline{z}$, and satisfies $\xi = \mathcal{O}(e^{\varepsilon t})~~as~t \to -\infty$.

Furthermore,  we obtain
\[
f'(\xi)
= -\lambda p \left( \Lambda + \xi \right)^{-(p+1)} + \lambda p \Lambda^{-(p+1)}
= \mathcal{O}(e^{\varepsilon t}).~~~~~ as ~~t \to -\infty
\]

 We first derive an estimate for $w$.
 \begin{lem}
Assume that $w(t,\theta)$ is a solution of \eqref{qcxx}. Then
\[
\max_{S^{N-1}} |w(t,\theta)| \leq 
\sigma_1^{(k)}, 
\]
\end{lem}
\begin{proof}

\text{For notational convenience, we fix }$k_0$ \text{such that } 
$\delta^{(k_0-1)} < \alpha < \delta^{(k_0)}$.

We expand $w$ as
\[
w(t, \theta)=\sum_{k=1}^{\infty} \sum_{j=0}^{m_k} w_j^k(t)\, Q_j^k(\theta).
\]

Since the first item of $z(t,\theta) - \overline{z}(t)$ vanishes, we have $m_0 = 1$ and $w_1^{(0)}(t)\equiv 0$. $w_j^{(k)}(t)$ satisfies the ODE
\begin{equation}
\left(w_j^{(k)}\right)'' 
+ \left( N-2 + 2\frac{\alpha+2}{p+1} \right) \left(w_j^{(k)}\right)'
+ \left[ (\alpha+2)\left( N-2 + \frac{\alpha+2}{p+1} \right) - \lambda_k \right] w_j^{(k)} 
= g_j^{(k)}(t),
\tag{12}\label{qdd}
\end{equation}
where
\[
\begin{aligned}
g_j^{(k)}(t)
&= \int_{\mathbb{S}^{N-1}} \left( f(z(t,\theta)) - \overline{f(z(t,\theta))} \right) Q_j^{(k)}(\theta)\, d\theta\\
&= \int_{\mathbb{S}^{N-1}} \left( f(z(t,\theta)) - f(\overline{z}(t)) \right) Q_j^{(k)}(\theta)\, d\theta
 - \int_{\mathbb{S}^{N-1}} \overline{ \left( f(z(t,\theta)) - f(\overline{z}(t)) \right) } Q_j^{(k)}(\theta)\, d\theta\\
&= \int_{\mathbb{S}^{N-1}} \left( f(z(t,\theta)) - f(\overline{z}(t)) \right) Q_j^{(k)}(\theta)\, d\theta.
\end{aligned}
\]

For $k \geq 1$, we observe that
\[
\| w \|_{L^2(\mathbb{S}^{N-1})}^2 =\sum_{k=1}^{\infty} \sum_{j=1}^{m_k} \left( w_j^{(k)}(t) \right)^2, \quad | f(z) - f(\overline{z})|=\sum_{k=1}^{\infty} \sum_{j=0}^{m_k}\left( g_j^{(k)}(t) \right)^2 .
\]

Since $f(z) - f(\overline{z}) = f'(\xi)\, w$ and $f'(\xi) = \mathcal{O}(e^{\varepsilon t})$, it follows that
\begin{equation}
\sum_{k=1}^{\infty} \sum_{j=1}^{m_k}\left( g_j^{(k)}(t) \right)^2
= \mathcal{O}(e^{2\varepsilon t}) \sum_{k=1}^{\infty} \sum_{j=1}^{m_k}\left( w_k^{(j)}(t) \right)^2.
\tag{13}\label{2222}
\end{equation}

On the other hand, applying \eqref{qdd}, we obtain for $T \ll -1$ and $t < T$,
\begin{equation}
\begin{aligned}
w_j^k(t)
= A_{j,1}^k e^{\sigma_1^{(k)} t}
+ A_{j,2}^k e^{\sigma_2^{(k)} t}
+ B_{j,1}^k \int_t^T e^{\sigma_1^{(k)}(t-s)} g_j^k(s)\, ds
- B_{j,2}^k \int_{-\infty}^t e^{\sigma_2^{(k)}(t-s)} g_j^k(s)\, ds,
\end{aligned}
\tag{14}\label{fg}
\end{equation}
where $\left| B_{j,1}^k \right|$ and $\left| B_{j,2}^k \right|$ are positive constants.

For $k \geq k_0$, we have $\sigma_1^{(k)} > 0$ and $\sigma_2^{(k)} < 0$.  
Since $w_j^k(t)\to 0$ as $t\to -\infty$, it follows that $A_{j,2}^k = 0$. Consequently,
\[
w_j^k(T)
= A_{j,1}^k e^{\sigma_1^{(k)} T}
- B_{j,2}^k \int_{-\infty}^T e^{\sigma_2^{(k)}(T-s)} g_j^k(s)\, ds,
\]
where $A_{j,1}^k = O\left(e^{-\sigma_1^{(k)} T}\right)$.

Therefore,
\begin{equation*}
\begin{aligned}
w_j^k(t)
= O\left(e^{\sigma_1^{(k)}(t-T)}\right)
+ B_{j,1}^k \int_t^T e^{\sigma_1^{(k)}(t-s)} g_j^k(s)\, ds
- B_{j,2}^k \int_{-\infty}^t e^{\sigma_2^{(k)}(t-s)} g_j^k(s)\, ds.
\end{aligned}
\end{equation*}

We may choose $\delta>0$ sufficiently small such that
\[
\begin{aligned}
\left[ w_j^k(t) \right]^2
\leq\; &
O\left(e^{2 \sigma_1^{(k)}(t-T)}\right)
+ 4(B_{j,1}^k)^2
\left( \int_t^T e^{\delta(t-s)}\, ds \right)
\left( \int_t^T e^{(2\sigma_1^{(k)} - \delta)(t-s)} (g_j^k(s))^2\, ds \right)
\\
&\quad
+ 4(B_{j,2}^k)^2
\left( \int_{-\infty}^t e^{-\delta(t-s)}\, ds \right)
\left( \int_{-\infty}^t e^{(2\sigma_2^{(k)} + \delta)(t-s)} (g_j^k(s))^2\, ds \right)
\\[1mm]
\leq\; &
C e^{2\sigma_1^{(k)}(t-T)}
+ C_\delta \int_t^T e^{(2\sigma_1^{(2)} - \delta)(t-s)} (g_j^k(s))^2\, ds
+ C_\delta \int_{-\infty}^t e^{(2\sigma_2^{(2)} + \delta)(t-s)} (g_j^k(s))^2\, ds,
\end{aligned}
\]
where $C > 0$ and $C_\delta > 0$.

For $k < k_0$,

\[
\begin{aligned}
w_j^k(t)= & A_{j, 1}^k e^{\sigma_1^{(k)} t}+A_{j, 2}^k e^{\sigma_2^{(k)} t} -B_{j, 1}^k \int_{-\infty}^t e^{\sigma_1^{(k)}(t-s)} g_j^k(s) d s-B_{j, 2}^k \int_{-\infty}^t e^{\sigma_2^{(k)}(t-s)} g_j^k(s) d s.
\end{aligned}
\]

In this case,
\[
\left|B_{j,1}^k\right|=\left|B_{j,2}^k\right|
= \left| \frac{1}{\sigma_2^{(k)} - \sigma_1^{(k)}} \right|.
\]

Since $\sigma_1^{(k)} < 0$ , $\sigma_2^{(k)} < 0$ and  $w_j^k(t)\to 0$ as $t\to -\infty$, we get that  
$A_{j,1}^k = A_{j,2}^k = 0$. Hence,
\begin{equation}
w_j^k(t)
= -B_{j,1}^k \int_{-\infty}^t e^{\sigma_1^{(k)}(t-s)} g_j^k(s)\, ds
 - B_{j,2}^k \int_{-\infty}^t e^{\sigma_2^{(k)}(t-s)} g_j^k(s)\, ds.
\tag{15}\label{13}
\end{equation}

Moreover,
\begin{equation}
\left( w_j^k(t) \right)^2 = O\left( e^{4\varepsilon t} \right).
\tag{16}\label{15}
\end{equation}

Observe that
\[
\left( g_j^k(t) \right)^2
\le C \, \| f'(\xi) w \|_{L^2(S^{N-1})}^2
\le C e^{4\varepsilon t}.
\]

Therefore,
\[
\begin{aligned}
\sum_{k=k_0}^{\infty} \sum_{j=1}^{m_k} \left( w_j^k(t) \right)^2
\leq &
\; C \sum_{k=k_0}^{\infty} \sum_{j=1}^{m_k} e^{2 \sigma_1^{(k)}(t-T)}
+ C_\delta \int_t^T e^{(2\sigma_1^{(k_0)} - \delta)(t-s)}
     \sum_{k=k_0}^{\infty} \sum_{j=1}^{m_k} \left( g_j^k(s) \right)^2\, ds
\\
&\; + C_\delta \int_{-\infty}^t e^{(2\sigma_2^{(k_0)} + \delta)(t-s)}
     \sum_{k=k_0}^{\infty} \sum_{j=1}^{m_k} \left( g_j^k(s) \right)^2\, ds
\\[1mm]
\leq &
\; C \sum_{k=k_0}^{\infty} \sum_{j=1}^{m_k} e^{2 \sigma_1^{(k)}(t-T)}
+ C \int_t^T e^{(2\sigma_1^{(k_0)} - \delta)(t-s)} e^{4\varepsilon s}\, ds
\\
&\; + C \int_t^T e^{(2\sigma_1^{(k_0)} - \delta)(t-s)}
     e^{2\varepsilon s} \sum_{k=k_0}^{\infty} \sum_{j=1}^{m_k}
     \left( w_j^k(s) \right)^2\, ds
\\
&\; + C \int_{-\infty}^t e^{(2\sigma_2^{(k_0)} + \delta)(t-s)} e^{4\varepsilon s}\, ds
\\
&\; + C \int_{-\infty}^t e^{(2\sigma_2^{(k_0)} + \delta)(t-s)}
     e^{2\varepsilon s} \sum_{k=k_0}^{\infty} \sum_{j=1}^{m_k}
     \left( w_j^k(s) \right)^2\, ds.
\end{aligned}
\]

Note that
\[
\sum_{k=k_0}^{\infty} \sum_{j=1}^{m_k} \left( g_j^k(t) \right)^2
\le \| f'(\xi) w \|_{L^2(S^{N-1})}^2
   - \sum_{k=1}^{k_0-1} \sum_{j=1}^{m_k} \left( g_j^k(t) \right)^2,
\]
and
\begin{equation}
\| f'(\xi) w \|_{L^2(S^{N-1})}^2
= O\left( e^{2\varepsilon t} \right)
  \| w \|_{L^2(S^{N-1})}^2
= O\left( e^{2\varepsilon t} \right)
  \left[ 
     \sum_{k=k_0}^{\infty} \sum_{j=1}^{m_k} \left( w_j^k(t) \right)^2
     + \sum_{k=1}^{k_0-1} \sum_{j=1}^{m_k}\left(w_j^k(t)\right)^2\right].
\end{equation}
Since
\[
\lim_{k\to\infty}
\frac{
m_{k+1} e^{2(\sigma_1^{(k+1)}-\sigma_1^{(k_0)})(t-T)}
}{
m_k e^{2(\sigma_1^{(k)}-\sigma_1^{(k_0)})(t-T)}
}
=
\lim_{k\to\infty}
\left[
\frac{m_{k+1}}{m_k}
e^{2(\sigma_1^{(k+1)}-\sigma_1^{(k)})(t-T)}
\right]
=
e^{(t-T)} < 1,
\]
we obtain
\[
\sum_{k=k_0}^{\infty} \sum_{j=1}^{m_k}
e^{2\sigma_1^{(k)}(t-T)}
=
\sum_{k=k_0}^{\infty} m_k e^{2\sigma_1^{(k)}(t-T)}
=
O\!\left( e^{2\sigma_1^{(k_0)}(t-T)} \right),
\]
Let
\[
[W(t)]^2
=
\sum_{k=k_0}^{\infty}
\sum_{j=1}^{m_k}
\left( w_j^k(t) \right)^2.
\]
Then
\[
\begin{aligned}
{[W(t)]}^2 \le &
\; C e^{2\sigma_1^{(k_0)}(t-T)}
+ C e^{4\varepsilon t}
+ C\int_t^T e^{(2\sigma_1^{(k_0)}-\delta)(t-s)} e^{4\varepsilon s}\, ds
\\
&\quad
+ C\int_t^T e^{(2\sigma_1^{(k_0)}-\delta)(t-s)}
       e^{2\varepsilon s} [W(s)]^2\, ds
+ C\int_{-\infty}^t e^{(2\sigma_2^{(k_0)}+\delta)(t-s)}
       e^{2\varepsilon s} [W(s)]^2\, ds.
\end{aligned}
\]

Next we consider two cases:  
(i) $4\varepsilon \ge 2\sigma_1^{(2)} - \delta$,  
(ii) $4\varepsilon < 2\sigma_1^{(2)} - \delta$.

For the first case, we first assume $4\varepsilon > 2\sigma_1^{(2)} - \delta$.  
We have
\begin{equation}
\begin{aligned}
{[W(t)]}^2 \le\;
& C e^{(2\sigma_1^{(2)}-\delta)(t-T)}
+ C\int_t^T e^{(2\sigma_1^{(2)}-\delta)(t-s)}
  e^{2\varepsilon s} [W(s)]^2\, ds
\\
&\quad
+ C\int_{-\infty}^t e^{(2\sigma_2^{(k_0)}+\delta)(t-s)}
  e^{2\varepsilon s} [W(s)]^2\, ds.
\end{aligned}
\tag{17}\label{16}
\end{equation}

Define
$$
K_1(t)
=
\int_t^T e^{(2\sigma_1^{(k_0)}-\delta)(t-s)} [W(s)]^2\, ds,
\qquad
K_2(t)
=
\int_{-\infty}^t e^{(2\sigma_2^{(k_0)}+\delta)(t-s)} [W(s)]^2\, ds.
$$

For $T$ sufficiently small, we compute that
$$
\begin{aligned}
\left( K_2 - K_1 \right)'(t)
&=
(2\sigma_2^{(k_0)}+\delta) K_2(t)
-
(2\sigma_1^{(k_0)}-\delta) K_1(t)
+ 2[W(t)]^2
\\
&\le
(2\sigma_2^{(k_0)}+\delta) K_2(t)
-
(2\sigma_1^{(k_0)}-\delta) K_1(t)
+ C e^{2\varepsilon T}\left( K_1(t)+K_2(t) \right)
+ C e^{(2\sigma_1^{(k_0)}-\delta)(t-T)}
\\
&\le C e^{(2\sigma_1^{(k_0)}-\delta)(t-T)}.
\end{aligned}
$$

Note  $\sigma_2^{(k_0)} < 0$ and $\sigma_1^{(k_0)} > 0$, and choose $\delta>0$ sufficiently small. Since $K_1(t) \to 0$ and $K_2(t) \to 0$ as $t\to -\infty$, we obtain for all $t < T$,
\begin{equation}
K_2(t) \le K_1(t) + C e^{(2\sigma_1^{(k_0)} - \delta)t}.
\tag{18}\label{17}
\end{equation}

Substituting \eqref{17} into \eqref{16}, we get at
\begin{equation}
[W(t)]^2
\le
C e^{(2\sigma_1^{(k_0)} - \delta)t}
+ C e^{2\varepsilon T} K_1(t).
\tag{19}\label{18}
\end{equation}

From \eqref{18}, it follows that
$$
K_1(t)
\le
\int_t^T
e^{(2\sigma_1^{(k_0)} - \delta)(t-s)}
\left[
C e^{(2\sigma_1^{(k_0)} - \delta)s}
+
C e^{2\varepsilon T} K_1(s)
\right] ds.
$$

Therefore,
\begin{equation*}
e^{-2(\sigma_1^{(k_0)} - \delta)t} \, K_1(t)
\le
C (T - t)
+
C e^{2\varepsilon T}
\int_t^T
K_1(s) \, e^{-(2\sigma_1^{(k_0)} - \delta)s}\, ds.
\end{equation*}

Let $F_1(t) = \int_t^T e^{-(2\sigma_1^{(k_0)} - \delta)s}\, ds$. Then
$$
- F_1'(t) \le C (T - t) + C e^{2\varepsilon T} F_1(t).
$$

Setting $\mu = C e^{2\varepsilon T}$, we obtain
$$
- (e^{\mu t} F_1(t))' \le C (T - t) e^{\mu t}.
$$

Integrating over $(t, T)$ yields
$$
F_1(t) \le \frac{C}{\mu^2} e^{-\mu(t - T)}.
$$

Hence,
$$
\begin{aligned}
K_1(t)
&\le
C e^{(2\sigma_1^{(k_0)} - \delta)t}(T - t)
+ C e^{2\mu} e^{(2\sigma_1^{(k_0)} - \delta)t} F_1(t)
\\
&\le
C e^{(2\sigma_1^{(k_0)} - \delta)t}(T - t)
+ \frac{C}{\mu}
   e^{(2\sigma_1^{(k_0)} - \delta - \mu)t} e^{\mu T}.
\end{aligned}
$$

Therefore,
$$
K_1(t) = O\left( e^{(2\sigma_1^{(k_0)} - \delta - \mu)t} \right).
$$

Consequently,
$$
[W(t)]^2
\le
C \mu \, e^{(2\sigma_1^{(k_0)} - \delta - \mu)t}.
$$

Meanwhile, by \eqref{15} and the assumption
$4\varepsilon \ge 2\sigma_1^{(k_0)} - \delta$, we also obtain
$$
\left( w_j^k(t) \right)^2
=
O\left( e^{(2\sigma_1^{(2)} - \delta)t} \right),
\qquad k = 1,2,\ldots, k_0 - 1.
$$

Hence, for all $t < T$,
$$
\sum_{k=1}^{\infty} \sum_{j=1}^{m_k}
\left( w_j^k(t) \right)^2
=
O\left( e^{(2\sigma_1^{(k_0)} - \delta - \mu)t} \right),
$$
 therefore

\begin{equation}
\| w \|_{L^2(S^{N-1})}
=
O\left(
   e^{\left( \sigma_1^{(k_0)} - \frac{\delta}{2}
             - \frac{\mu}{2} \right)t}
  \right).
\tag{20}\label{19}
\end{equation}

\begin{equation}
\max_{S^{N-1}} |w(t,\theta)|
\le
M e^{\left( \sigma_1^{(k_0)} - \frac{\delta}{2} - \frac{\mu}{2} \right)t},
\qquad \text{for all } t \in (-\infty, -1].
\tag{21}\label{20}
\end{equation}

We establish \eqref{20} only for $t \in (-\infty, T_*)$, ($T_* \le T$).  
The remaining part can be proved directly from the continuity of $w$.  
Define
$$
h(r,\theta) = w(t,\theta), \qquad r = e^t.
$$

Then $h(r,\theta)$ satisfies the equation
\begin{equation}
\Delta h
+ \frac{b_1 \, x\cdot \nabla h}{r^2}
+ \frac{b_2 h}{r^2}
- \frac{ f'(\xi) h - \overline{ f'(\xi) h } }{r^2}
= 0,
\qquad \text{in } B_{R_*} \setminus \{0\},
\tag{22}\label{21}
\end{equation}
where
$$
b_1 = 2\frac{\alpha+2}{p+1}, \qquad
b_2 = 2\left( N - 2 + \frac{\alpha+2}{p+1} \right),
\qquad R_* = e^{T_*}.
$$

For any $x_0 \in B_{R_*}\setminus\{0\}$, denote $r_0 = |x_0|$ and set  
$\Omega = B_{r_0/2}(x_0)$.  
We regard equation \eqref{21} as the linear equation appearing in Lemma 5.1 of reference \cite{3z}, where
$$
k = k_1 = 1, \qquad
h(x) \equiv 0, \qquad
|b| = \frac{4 b_1^2}{r_0^2}, \qquad
|c| = \frac{Q}{r_0^2},
$$
Here $Q = Q(h) > 0$.  
Thus, by applying Lemma 5.1 in reference \cite{3z} with $k_2 = Q / r_0^2$, and combining \eqref{19} with an argument similar to the proof of Theorem 5.1 in reference \cite{3z}, we conclude that there exists a positive constant
$$
M = M\!\left( \frac{k_1}{k},\, k_2 r_0^2 \right)
    = M(Q)
    = M(h),
$$
independent of $r_0$, such that
$$
\sup_{x \in B_{r_0/4}(x_0)}
|h(x)|
\le
M\, r_0^{\sigma_1^{(k_0)} - \frac{\delta}{2} - \frac{\mu}{2}}.
$$

In particular,
$$
|h(x_0)|
\le
M\, r_0^{\sigma_1^{(k_0)} - \frac{\delta}{2} - \frac{\mu}{2}},
\qquad
\max_{|x| = r} |h(x)|
\le
M\, r^{\sigma_1^{(k_0)} - \frac{\delta}{2} - \frac{\mu}{2}}.
$$

We obtain from estimate \eqref{20} that
\begin{equation}
\sum_{k=1}^{\infty} \sum_{j=1}^{m_k}
\left( g_j^k(t) \right)^2
=
O\!\left(
      e^{(2\sigma_1^{(k_0)} + 2\varepsilon - \delta - \mu)t}
   \right).
\tag{23}\label{22}
\end{equation}

By \eqref{13} and \eqref{22}, we know
\begin{equation}
w_j^k(t)
=
O\!\left(
     e^{(\sigma_1^{(k_0)} + \varepsilon - \frac{\delta}{2} - \frac{\mu}{2})t}
   \right)
=
O\!\left( e^{\sigma_1^{(k_0)} t} \right),
\qquad k \le k_0,
\tag{24}\label{23}
\end{equation}
since we may choose $\delta$ sufficiently small and $T$ sufficiently large so that  
$\delta < 2\varepsilon - \mu$.

Moreover, by using \eqref{22} we have
\[
\begin{aligned}
\sum_{j=1}^{m_{k_0}} |w_j^{k_0}(t)|
\le\;&
C e^{\sigma_1^{(k_0)} t}
+ C \int_t^T
 e^{\sigma_1^{(k_0)}(t-s)}
 \sum_{j=1}^{m_{k_0}} |g_j^{k_0}(s)|\, ds
\\
&\quad
+ C \int_{-\infty}^t
  e^{\sigma_2^{(k_0)}(t-s)}
  \sum_{j=1}^{m_{k_0}} |g_j^{k_0}(s)|\, ds
\\[2mm]
\le\;&
C e^{\sigma_1^{(k_0)} t}
+ C \int_t^T
    e^{\sigma_1^{(k_0)}(t-s)}
    e^{(\sigma_1^{(k_0)} + \varepsilon - \frac{\delta}{2} - \frac{\mu}{2})s}\, ds
\\
&\quad
+ C \int_{-\infty}^t
    e^{\sigma_2^{(k_0)}(t-s)}
    e^{(\sigma_1^{(k_0)} + \varepsilon - \frac{\delta}{2} - \frac{\mu}{2})s}\, ds
\\[2mm]
\le\;&
C e^{\sigma_1^{(k_0)} t}.
\end{aligned}
\]

Since $\delta < 2\varepsilon - \mu$, it follows for $t < T$ that
\begin{equation}
\sum_{j=1}^{m_{k_0}}
|w_j^{k_0}(t)|^2
\le
\left[
   \sum_{j=1}^{m_{k_0}} |w_j^{k_0}(t)|
\right]^2
\le
C e^{2\sigma_1^{(k_0)} t}.
\tag{25}\label{24}
\end{equation}

Similarly, for $k > k_0$ we obtain
\begin{equation*}
\begin{aligned}
\sum_{j=1}^{m_{k}}\left|w_j^{k}(t)\right| \leq&  C e^{\sigma_1^{(k)} t}+C \int_t^T e^{\sigma_1^{(k)}(t-s)} \sum_{j=1}^{m_{k}}\left|g_j^2(s)\right| d s 
 +C \int_{-\infty}^t e^{\sigma_2^{(k)}(t-s)} \sum_{j=1}^{m_{k}}\left|g_j^{k}(s)\right| d s \\ 
 \leq & C e^{\sigma_1^{(k)} t}+C \int_t^T e^{\sigma_1^{(k)}(t-s)} e^{\left(\sigma_1^{(k)}+\varepsilon-\frac{\delta}{2}-\frac{\mu}{2}\right) s} d s 
 +C \int_{-\infty}^t e^{\sigma_2^{(k)}(t-s)} e^{\left(\sigma_1^{(k)}+\varepsilon-\frac{\delta}{2}-\frac{\mu}{2}\right) s} d s \\
\leq & C e^{\sigma_1^{(k)} t} .
\end{aligned}
\end{equation*}

Hence
\begin{equation}
\sum_{k = k_0 + 1}^{\infty} \sum_{j=1}^{m_k}
\left( w_j^k(t) \right)^2
=
O\!\left( e^{\sigma_1^{(k_0)} t} \right).
\tag{26}\label{25}
\end{equation}

Since $\delta$ can be chosen sufficiently small, combining \eqref{23}--\eqref{25} yields
\begin{equation}
\sum_{k = 1}^{\infty} \sum_{j=1}^{m_k}
\left( w_j^k(t) \right)^2
=
O\!\left( e^{\sigma_1^{(k_0)} t} \right).
\tag{27}\label{26}
\end{equation}

The case $4\varepsilon = 2\sigma_1^{(2)} - \delta$ can be down similarly.
Indeed, by enlarging $\delta$ slightly to a number $\delta' > \delta$ such that  
$4\varepsilon > 2\sigma_1^{(2)} - \delta'$, and repeating the proof above with $\delta'$ in stead of $\delta$, we conclude \eqref{26}.

For case (ii), we have
$$
\begin{aligned}
{[W(t)]}^2
\le\;&
C e^{4\varepsilon t}
+ C \int_t^T
e^{(2\sigma_1^{(2)} - \delta)(t-s)}
e^{2\varepsilon s} [W(s)]^2\, ds
\\
&\quad
+ C \int_{-\infty}^t
e^{(2\sigma_2^{(2)} + \delta)(t-s)}
e^{2\varepsilon s} [W(s)]^2\, ds.
\end{aligned}
$$

Since $[W(t)]^2 = O(e^{2\varepsilon t})$, it follows that
$$
[W(t)]^2 \le C e^{4\varepsilon t},
\qquad t < T.
$$

Together with \eqref{15}, this also implies
$$
\| w \|_{L^2(S^{N-1})}
=
O\!\left( e^{2\varepsilon t} \right).
$$

Thus, by an argument similar to the proof of \eqref{20}, there exists a constant  
$M = M(w) > 0$ such that
\begin{equation}
\max_{S^{N-1}} |w(t,\theta)|
\le
M e^{2\varepsilon t},
\qquad t \in (-\infty, -1].
\tag{28}\label{27}
\end{equation}

From \eqref{13} and \eqref{27}, we obtain, for $k < k_0$,
$$
w_j^k(t)
=
O\!\left( e^{3\varepsilon t} \right).
$$

Therefore, we get the inequality
\begin{equation}
\begin{gathered}
{[W(t)]}^2
\le
C e^{2\sigma_1^{(k_0)}(t - T)}
+ C \int_t^T e^{(2\sigma_1^{(k_0)} - \delta)(t-s)} e^{6\varepsilon s}\, ds
+ C \int_t^T e^{(2\sigma_1^{(k_0)} - \delta)(t-s)}
     e^{2\varepsilon s} [W(s)]^2\, ds
\\
\quad
+ C \int_{-\infty}^t e^{(2\sigma_2^{(k_0)} + \delta)(t-s)} e^{6\varepsilon s}\, ds
+ C \int_{-\infty}^t e^{(2\sigma_2^{(k_0)} + \delta)(t-s)}
     e^{2\varepsilon s} [W(s)]^2\, ds.
\end{gathered}
\tag{29}\label{28}
\end{equation}

Note that:

$$
\sum_{k=k_0}^{\infty} \sum_{j=1}^{m_k}\left(g_j^k(t)\right)^2 \leq C e^{2 \epsilon t}\left([W(t)]^2+e^{6 \epsilon t}\right)
$$

We still consider two cases:  
(a) $6\varepsilon \ge 2\sigma_1^{(k_0)} - \delta$,  
(b) $6\varepsilon < 2\sigma_1^{(k_0)} - \delta$.

For case (a), using inequality \eqref{28} and an argument similar in case (i), we conclude that \eqref{26} holds.

For case (b), inequality \eqref{28} implies that
$$
[W(t)]^2 \le C e^{6\varepsilon t},
\qquad \text{for } t < T.
$$

Hence,
$$
\begin{gathered}
{[W(t)]}^2
\le
C e^{2\sigma_1^{(k_0)}(t - T)}
+ C \int_t^T e^{(2\sigma_1^{(k_0)} - \delta)(t-s)} e^{8\varepsilon s}\, ds
+ C \int_t^T e^{(2\sigma_1^{(k_0)} - \delta)(t-s)}
        e^{2\varepsilon s}[W(s)]^2\, ds
\\
\quad
+ C \int_{-\infty}^t e^{(2\sigma_2^{(k_0)} + \delta)(t-s)} e^{8\varepsilon s}\, ds
+ C \int_{-\infty}^t e^{(2\sigma_2^{(k_0)} + \delta)(t-s)}
        e^{2\varepsilon s}[W(s)]^2\, ds.
\end{gathered}
$$

We still consider into two cases:  
$8\varepsilon \ge 2\sigma_1^{(k_0)} - \delta$ and $8\varepsilon < 2\sigma_1^{(k_0)} - \delta$.  
Repeating the same steps as before, we eventually deduce that \eqref{26} holds.

Clearly, \eqref{26} implies
$$
\| w \|_{L^2(S^{N-1})}
=
O\!\left( e^{\sigma_1^{(k_0)} t} \right).
$$

Using the  argument similar as in Theorem 5.1 of reference \cite{3z}, we also obtain
$$
\max_{S^{N-1}} |w(t,\theta)|
\le
C e^{\sigma_1^{(k_0)} t},
\qquad t \in (-\infty,-1].
$$

This completes the proof of the lemma.
\end{proof}
\vspace{1em}

\begin{lem}
Let $\bar{z}(t)$ be defined as above. Then:

\[
|\bar{z}(t)| \leq
C e^{\min\left\{ 2\sigma_1^{(k)},\, \frac{2p - \alpha}{p+1} \right\} t}
\]
for $F \neq 0$.
\[
|\bar{z}(t)| \le
C e^{2\sigma_1^{(k)} t}
\]
for $F = 0$, $C > 0$ is a constant independent of $t$.
\end{lem}
\begin{proof} We only prove the case $F\neq0$,  
the other cases are similar and easier.  
 we only need to prove the estimate for $t$ sufficiently close to $-\infty$ and the other part can be obtained easily by the continuity of $\bar{z}(t)$.

We assume
$
\frac{2p - \alpha}{p+1} \le 2\sigma_1^{(k)},
$
while the opposite case
$
\frac{2p - \alpha}{p+1} > 2\sigma_1^{(k)}
$
can be treated similarly.

We know that $\bar{z}(t)$ satisfies the ODE
\[
\left\{
\begin{array}{l}
\bar{z}_{tt}
+
\left( N - 2 + 2\frac{\alpha+2}{p+1} \right) \bar{z}_t
+
(\alpha+2)\left( N - 2 + \frac{\alpha+2}{p+1} \right) \bar{z}
=
\overline{f(z)} + F e^{\frac{2p - \alpha}{p+1} t},
\quad t \in (-\infty, 0),
\\[1mm]
\bar{z}(t) = O(e^{\varepsilon t}) \to 0
\quad \text{as}~t \to -\infty,
\end{array}
\right.
\]
where $\varepsilon > 0$ as in Lemma 3.2.

Observe that
\[
\overline{f(z)} + F e^{\frac{2p - \alpha}{p+1} t}
=
O\big( \overline{(\bar{z} + w)^2} \big)
+ F e^{\frac{2p - \alpha}{p+1} t}
=
O\big( \overline{ \bar{z}^2 + 2\bar{z}w + w^2 } \big)
+ F e^{\frac{2p - \alpha}{p+1} t}.
\]

By Lemma 3.3
\[
\begin{aligned}
\overline{f(z)} + F e^{\frac{2p - \alpha}{p+1} t}
&=
O\left( \bar{z}^2
+ 2\bar{z}\, O(e^{\sigma_1^{(k)} t})
+ O(e^{2\sigma_1^{(k)} t}) \right)
+ F e^{\frac{2p - \alpha}{p+1} t}
\\
&=
O(e^{2\varepsilon t})
+ O(e^{(\sigma_1^{(k)}+\varepsilon)t})
+ O\!\left( e^{\min\{2\sigma_1^{(k)},\,\frac{2p - \alpha}{p+1}\} t} \right)
\\
&=
O(e^{2\varepsilon t})
+ O\left( e^{\frac{2p - \alpha}{p+1} t} \right).
\end{aligned}
\]
as $t\to -\infty$.

Since  
\[
\frac{2p - \alpha}{p+1} \le 2\sigma_1^{(k)}, 
\qquad
0 < \varepsilon \le \sigma_1^{(k)},
\]
we consider two cases:  
(i) $2\varepsilon \ge \frac{2p - \alpha}{p+1}$,  
(ii) $2\varepsilon < \frac{2p - \alpha}{p+1}$.

 Case (i):  \( 2\varepsilon \ge \frac{2p - \alpha}{p+1} \)

In this case,
\[
\overline{f(z)} + F e^{\frac{2p - \alpha}{p+1} t}
=
O\left( e^{\frac{2p - \alpha}{p+1} t} \right).
\]

Using the representation formula for solutions of the ODE for $\bar{z}$,  
and following the similar argument as in the proof of Theorem 2.1, we obtain
\[
\bar{z}(t)
=
O\left( e^{\frac{2p - \alpha}{p+1} t} \right).
\]
This completes the proof of Lemma 3.4 in this case.

 Case (ii):  \( 2\varepsilon < \frac{2p - \alpha}{p+1} \)

We have
\[
\overline{f(z)} + F e^{\frac{2p - \alpha}{p+1} t}
=
O(e^{2\varepsilon t}).
\]

Using again the ODE theory for $\bar{z}$ and an argument similar to the proof of Theorem 2.1, we conclude that
\[
\bar{z}(t) = O(e^{2\varepsilon t}).
\]

By Lemma 3.3, we have
$$
z(t,\theta)
= w(t,\theta) + \bar{z}(t)
= O\!\left( e^{\sigma_1^{(k)} t} \right)
  + O\!\left( e^{2\varepsilon t} \right).
$$

 we obtain
\[
z(t,\theta) = O\!\left( e^{\sigma_1^{(k)} t} \right),
\qquad
f(z) = O(z^2) = O\!\left( e^{2\sigma_1^{(k)} t} \right).
\]
When $2\sigma_1^{(k)} \le \varepsilon$.

Thus
\[
\overline{f(z)} + F e^{\frac{2p-\alpha}{p+1} t}
=
O\!\left( e^{\frac{2p-\alpha}{p+1} t} \right).
\]

Using the ODE for $\bar{z}$, we obtain
\[
\bar{z}(t)
= O\!\left( e^{\frac{2p-\alpha}{p+1} t} \right).
\]

Also we obtain
\[
z(t,\theta) = O\!\left( e^{2\varepsilon t} \right),
\]
When $2\varepsilon < \sigma_1^{(k)}$,
and therefore
\[
\overline{f(z)} + F e^{\frac{2p-\alpha}{p+1} t}
=
O\!\left( e^{4\varepsilon t} \right)
+ O\!\left( e^{\frac{2p-\alpha}{p+1} t} \right).
\]

We still consider two cases:  
(a) $4\varepsilon \ge \frac{2p-\alpha}{p+1}$,  
(b) $4\varepsilon < \frac{2p-\alpha}{p+1}$.  
Repeating the similar steps as above then we yield the desired estimate.
\end{proof} 
We now prove \textbf{Lemma 3.2}.  
From the arguments given above, observe that  
when $\alpha < (\ge) \, 2p - (p+1)\sigma_1^{(k)}$,  
we have $\sigma_1^{(k)} < (\ge)\, \frac{2p - \alpha}{p+1}$.  
Since $z(t,\theta) = w(t,\theta) + \bar{z}(t)$, the conclusion follows immediately.

\vspace{1em}

Finally, we proceed to prove \textbf{Theorem 3.1}.
\begin{proof}
(i) For the case $F = 0$ and $\delta^{(k-1)} < \alpha < \delta^{(k)}$, we set
\begin{equation}
\rho_1 = \sigma_1^{(k)},\;
\rho_2 = \sigma_1^{(k+1)},\; \ldots,\;
\rho_m = \sigma_1^{(k+m-1)},\; \ldots
\tag{30}\label{30}
\end{equation}

(ii) For the case $F \neq 0$ and $\delta^{(k-1)} < \alpha < \delta^{(k)}$,  
if there exists some $j \ge k$ such that
\[
\sigma_1^{(j)} < \frac{2p - \alpha}{p+1} < \sigma_1^{(j+1)},
\]
we set
\begin{equation}
\rho_1 = \sigma_1^{(k)},\;
\rho_2 = \sigma_1^{(k+1)},\; \ldots,\;
\rho_{j-k+1} = \sigma_1^{(j)},\;
\rho_{j-k+2} = \frac{2p - \alpha}{p+1},\;
\rho_{j-k+3} = \sigma_1^{(j+1)},\; \ldots
\tag{31}\label{34}
\end{equation}

~~~~For the case $F \neq 0$ and $\delta^{(k-1)} < \alpha < \delta^{(k)}$,  
if there exists some $j \ge k$ such that  
\[
\sigma_1^{(j)} < \frac{2p - \alpha}{p+1} = \sigma_1^{(j+1)}
\quad \text{or} \quad
\sigma_1^{(j)} = \frac{2p - \alpha}{p+1} < \sigma_1^{(j+1)},
\]
we set
\begin{equation}
\rho_1 = \sigma_1^{(k)},\;
\rho_2 = \sigma_1^{(k+1)},\; \ldots,\;
\rho_{j-k+1} = \sigma_1^{(j)},\;
\rho_{j-k+2} = \sigma_1^{(j+1)},\; \ldots
\tag{32}\label{35}
\end{equation}

~~~~~For the case $F \neq 0$ and $\delta^{(k-1)} < \alpha < \delta^{(k)}$,  
if
\[
\frac{2p - \alpha}{p+1} < \sigma_1^{(k)},
\]
we set
\begin{equation}
\rho_1 = \frac{2p - \alpha}{p+1},
\quad
\rho_2 = \sigma_1^{(k)},\; \ldots,\;
\rho_m = \sigma_1^{(m - 2 + k)},\; \ldots
\tag{33}\label{36}
\end{equation}

The sequence $\{\rho_i\}_{i \ge 1}$ is strictly increasing and diverges to $+\infty$ under the above assumptions,

We prove the subcase that $\delta^{(0)}<\alpha<\delta^{(1)}$ of case \eqref{30} at first, and then indicate the proof of the remaining cases.

\vspace{1em}

We start from the identity
$$
\mathcal{L} z = f(z),
$$
where
$$
\mathcal{L} z
=
z_{tt}
+ \left( N - 2 + 2\frac{\alpha+2}{p+1} \right) z_t
+ \Delta_{S^{N-1}} z
+ (\alpha+2)\left( N - 2 + \frac{\alpha+2}{p+1} \right) z,
$$
and
$$
f(z)
=
\lambda (\Lambda + z)^{-p}
- \lambda \Lambda^{-p}
+ \lambda p z \Lambda^{-(p+1)}.
$$

For $-\Delta_\theta Q_i = \lambda_i Q_i$ $(i \ge 0)$,  
and for clarity of presentation, we write the eigenvalues with multiplicities:
$$
\lambda_0 = 0,\qquad
\lambda_1 = \cdots = \lambda_n = 1,\qquad
\lambda_{n+1} = 2n,\quad \ldots
$$
We fix $\{Q_i\}$ as an orthonormal basis of $L^2(\mathbb{S}^{n-1})$.

For each fixed $i \ge 0$ and any twice differentiable function $\psi = \psi(t)$, we define
$$
\mathcal{L}(\psi Q_i) = (L_i \psi)\, Q_i.
$$
Since $-\Delta_\theta Q_i = \lambda_i Q_i$, we obtain
$$
L_i \psi
=
\psi_{tt}
+ \left( N - 2 + 2\frac{\alpha+2}{p+1} \right)\psi_t
- \lambda_i \psi
+ (\alpha+2)\left( N - 2 + \frac{\alpha+2}{p+1} \right)\psi.
$$
\begin{lem}
There exist two sequences $\{\rho_i\}_{i\ge 1}$ and $\{\tau_i\}_{i\ge 1}$,  
tending respectively to $+\infty$ and $-\infty$,  
such that for every $i \ge 1$,  
$\operatorname{Ker}(L_i)$ has a basis $\{e^{\rho_i t},\, e^{\tau_i t}\}$.
\end{lem}
We now introduce the notion of an index set.  
Let $\{\rho_i\}_{i\ge 1}$ denote the sequence appearing in Lemma 3.5
 (i.e. the sequence $\{\rho_i\}$ in \eqref{30} for $k=1$ with multiplicities considered),  
which is strictly increasing and diverges to $+\infty$.
\vspace{1em}

We define the index set
$$
\mathcal{I}
=
\left\{
\sum_{i\ge 1} m_i \rho_i
\;\middle|\;
m_i \text{ are positive integers with only finitely many } m_i > 0
\right\}.
$$

In other words, $\mathcal{I}$ consists of all finite positive-integer linear combinations of the $\rho_i$.  
It's possible that a given $\rho_i$ may itself be representable as a positive-integer linear combination of  
$\rho_1,\ldots,\rho_{i-1}$.

We now give another lemma.
\begin{lem}
If $Q_k$ and $Q_l$ are spherical harmonics of degrees $k$ and $l$ respectively, then
$$
Y_k Y_l = \sum_{i=0}^{k+l} Z_i,
$$
where each $Z_i$ is a spherical harmonic of degree $i$ $(i = 0,1,\ldots,k+l)$.
\end{lem}
\emph{Proof.}
We use polar coordinates $(r,\theta)$ on $\mathbb{R}^n$.  
Then
$u_k(x) = r^k Q_k(\theta)$ and $u_l(x) = r^l Q_l(\theta)$
are homogeneous harmonic polynomials of degrees $k$ and $l$ respectively.  
Hence $u_k u_l$ is a homogeneous polynomial of degree $k+l$.  

By the decomposition theorem for homogeneous polynomials stated in reference \cite{5z}, we have
$$
u_k(x) u_l(x)
=
v_{k+l}(x) + |x|^2 v_{k+l-2}(x)
+ \cdots
+ |x|^{\,k+l-\tau} v_\tau(x),
$$
where $\tau = 1$ if $k+l$ is odd and $\tau = 0$ if $k+l$ is even,  
and each $v_i$ is a homogeneous harmonic polynomial of degree $i$  
$(i = k+l,\, k+l-2,\,\ldots,\,\tau)$.  
Restricting the above identity to the unit sphere yields Lemma 3.6.

We recall that $\{Q_i\}$ are the eigenfunctions of $-\Delta_\theta$ and it forms an orthonormal eigenbasis in
$L^2(\mathbb{S}^{n-1})$.
The corresponding eigenvalues $\{\lambda_i\}$ are increasing.  
Thus each $Q_i$ is a spherical harmonic of degree $\deg(Q_i)$,  
and  we have $\deg(Q_i) \le \deg(Q_j)$ for $i \le j$ . 
Here $Q_0$ is constant and $Q_1,\ldots,Q_n$ are degree-1 spherical harmonics.  
Since $z(t,\theta) = O(e^{\varepsilon t})$ and$f(z) = \sum_{i=2}^{\infty} c_i z^i$. We know
\[
|\mathcal{L}z| = |f(z)| \le C z^2,
\]
we now decompose the index set $\mathcal{I}$.

Define
\begin{equation}
\mathcal{I}_\rho = \{\rho_i : i \ge 1\},
\tag{34}\label{36.1}
\end{equation}
and
\begin{equation}
\mathcal{I}_{\widetilde{\rho}}
=
\left\{
\sum_{i=1}^r n_i \rho_i :
n_i \in \mathbb{Z}_+,\;
\sum_{i=1}^r n_i \ge 2
\right\}.
\tag{35}\label{36.2}
\end{equation}

We assume that the sequence in $\mathcal{I}_{\widetilde{\rho}}$ is
$\{\widetilde{\rho}_i\}_{i \ge 1}$, which is strictly increasing with
$\widetilde{\rho}_1 = 2\rho_1$.  
We first consider the case where
$\mathcal{I}_\rho \cap \mathcal{I}_{\widetilde{\rho}} = \emptyset$;  
that is, no $\rho_i$ can be expressed as a positive-integer linear combination of
$\rho_1,\ldots,\rho_{i-1}$ (except for the trivial identity
$\rho_{\widetilde{i}} = \rho_i$).  
In this situation, the elements of $\mathcal{I}$ may be arranged as follows:
\begin{equation}
\rho_1 \le \cdots \le \rho_{r_1}
< \widetilde{\rho}_1
< \cdots < \widetilde{\rho}_{l_1}
< \rho_{r_1+1}
\le \cdots \le \rho_{r_2}
< \widetilde{\rho}_{l_1+1}
< \cdots .
\tag{36}\label{28q}
\end{equation}

For each $\widetilde{\rho}_i \in \mathcal{I}_{\widetilde{\rho}}$,  
we consider nonnegative integers $n_1,\ldots,n_{r_1}$ such that
\begin{equation}
n_1 + \cdots + n_{r_1} \ge 2,
\qquad
n_1 \rho_1 + \cdots + n_{r_1} \rho_{r_1}
= \widetilde{\rho}_i .
\tag{37}\label{29q}
\end{equation}
Clearly, only finitely many such $(n_1,\ldots,n_{r_1})$ exist.

Define
$$
\begin{aligned}
\widetilde{R}_i
=
\max\Big\{
& n_1 \deg(Q_1)
+ n_2 \deg(Q_2)
+ \cdots
+ n_{r_1} \deg(Q_{r_1})
: 
\\[1mm]
& \left.
(n_1,\ldots,n_{r_1})
\text{ are nonnegative integers satisfying \eqref{29q}}
\right\},
\end{aligned}
$$
and
\begin{equation}
\widetilde{M}_i
=
\max\{\, m : \deg(Q_m) \le \widetilde{R}_i \,\}.
\tag{38}\label{30q}
\end{equation}
We obtain From Lemma 3.2 that
$$
z=O\left(e^{\rho_1 t}\right)
$$
Hence
$$
|\mathcal{L}(z)| = O(e^{2\rho_1 t}) = O(e^{\widetilde{\rho}_1 t}).
$$

We now proceed in several steps to establish the case
$\mathcal{I}_\rho \cap \mathcal{I}_{\widetilde{\rho}} = \varnothing$.

\textbf{Step 1.}
Observe that $\rho_{r_1} < \widetilde{\rho}_1 = 2\rho_1$.  
Thus, by Lemma A.8 in \cite{4z}
(although the statement there is for $t \to +\infty$, the conclusion clearly remains valid as $t \to -\infty$; we will not repeat this remark later),  
there exists a function $\eta_1$ such that
$$
z = \eta_1 + O(e^{\widetilde{\rho}_1 t}),
$$
where
\begin{equation}
\eta_1(t,\theta)
=
\sum_{i=1}^{r_1} c_i Q_i(\theta) e^{\rho_i t}.
\tag{39}\label{31q}
\end{equation}

Define
\begin{equation}
z_1 = z - \eta_1.
\tag{40}\label{32q}
\end{equation}
Then $\mathcal{L}\eta_1 = 0$, $\mathcal{L}z_1 = f(z)$, and
\begin{equation}
z_1 = O(e^{\widetilde{\rho}_1 t}).
\tag{41}\label{33q}
\end{equation}

\textbf{Step 2.}
We next show that there exists a function $\widetilde{\eta}_1$ such that
\begin{equation}
\widetilde{z}_1
=
z_1 - \widetilde{\eta}_1
=
z - \eta_1 - \widetilde{\eta}_1,
\tag{42}\label{34q}
\end{equation}
and
\begin{equation}
\mathcal{L}\widetilde{z}_1
=
O(e^{\widetilde{\rho}_{l_1+1} t}).
\tag{43}\label{35q}
\end{equation}

We claim that $\widetilde{\eta}_1$ has the structure
\begin{equation}
\widetilde{\eta}_1(t,\theta)
=
\sum_{i=1}^{l_1}
\left\{
\sum_{m=0}^{\widetilde{M}_i} c_{im} Q_m(\theta)
\right\}
e^{\widetilde{\rho}_i t},
\tag{44}\label{36q}
\end{equation}
where $\widetilde{M}_i$ is as defined in \eqref{30q},  
and $c_{im}$ are constants.  
This provides the next level of expansion.

To prove this, fix $\widetilde{\eta}_1$ as above and define $\widetilde{z}_1$.  
Then
\begin{equation}
\mathcal{L}\widetilde{z}_1
=
f(z) - \mathcal{L}\widetilde{\eta}_1.
\tag{45}\label{37q}
\end{equation}

Note that $3\rho_1 \in \mathcal{I}_{\widetilde{\rho}}$.  
We divide the discussion into two cases.

\textbf{Case 1.}
Assume $\rho_{r_1+1} < 3\rho_1$.  
Then $\widetilde{\rho}_{l_1} < \rho_{r_1+1} < 3\rho_1$,  
which implies $\widetilde{\rho}_{l_1+1} \le 3\rho_1$.

Observe that
$$
f(z)
=
f(z_1 + \eta_1)
=
\sum_{i=2}^{\infty} c_i (z_1 + \eta_1)^i.
$$

we get from \eqref{32q} that 
$$
z_1^2 \le C e^{4\rho_1 t},
\qquad
|z_1 \eta_1| \le C e^{3\rho_1 t}.
$$

From the expression of $\eta_1$ in \eqref{31q}, we obtain
$$
\sum_{i=2}^{\infty} c_i \eta_1^i
=
\sum_{n_1+\cdots+n_{r_1}\ge 2}
a_{n_1\cdots n_{r_1}}
e^{(n_1\rho_1 + \cdots + n_{r_1}\rho_{r_1}) t}
Q_1^{n_1}\cdots Q_{r_1}^{n_{r_1}}.
$$

By the definition of $\mathcal{I}_{\widetilde{\rho}}$,  
there exists that $\widetilde{\rho}_i = n_1\rho_1 + \cdots + n_{r_1}\rho_{r_1}$.  
Thus, by Lemma 3.6,
$$
\sum_{i=2}^{\infty} c_i \eta_1^i=\sum_{i=1}^{\infty}\left\{\sum_{m=0}^{\widetilde{M}_i} a_{i m}Q_m(\theta)\right\} e^{\widetilde{\rho}_i t}
$$
We now truncate the summation on the right-hand side at the finite index $l_1$ and denote the expression by $I_1$. Then
$$
I_1
=
\sum_{i=1}^{l_1}
\left\{
\sum_{m=0}^{\widetilde{M}_i}
a_{im} Q_m(\theta)
\right\}
e^{\widetilde{\rho}_i t}.
$$
Hence,
$$
f(z) = I_1 + O(e^{\widetilde{\rho}_{l_1+1} t}).
$$

Therefore, by \eqref{37q},
$$
\mathcal{L} z_1
=
\mathcal{L} \widetilde{\eta}_1 - I_1
+ O(e^{\widetilde{\rho}_{l_1+1} t}).
$$

We choose $\widetilde{\eta}_1$ to be the form 
$$
\widetilde{\eta}_1(t,\theta)
=
\sum_{i=1}^{l_1}
\sum_{m=0}^{\widetilde{M}_i}
\widetilde{\eta}_{im}(t) Q_m(\theta).
$$

To solve the equation $\mathcal{L}\widetilde{\eta}_1 = I_1$. We impose
\begin{equation}
L_m \widetilde{\eta}_{im}
=
a_{im} e^{-\widetilde{\rho}_i t}.
\tag{46}\label{38q}
\end{equation}
for each $1 \le i \le l_1$ and $0 \le m \le \widetilde{M}_i$.

Since $\rho_m \ne \widetilde{\rho}_i$ for every $m \neq i$, we can get from
Lemma A.2 and Remark A.5 in \cite{4z}(constants are viewed as special periodic functions) that
\begin{equation}
\widetilde{\eta}_{im}(t)
=
c_{im} e^{\widetilde{\rho}_i t}.
\tag{47}\label{39q}
\end{equation}

Thus we have obtained explicit formulas for $\widetilde{\eta}_1$ and $z_1$.  
Moreover, by \eqref{33q} and \eqref{36q}, we obtain
$$
\widetilde{z}_1 = O(e^{\widetilde{\rho}_1 t}).
$$

\textbf{Case 2.}
Assume now that $\rho_{r_1+1} > 3\rho_1$.  
Then $\widetilde{\rho}_{l_1} \ge 3\rho_1$.  
Let $n_1$ be the largest integer satisfying $\widetilde{\rho}_{n_1} < 3\rho_1$.  
Then $\widetilde{\rho}_{n_1+1} = 3\rho_1$.

We repeat the argument of Case 1 with $n_1$ replacing $l_1$,  
and redefine $I_1$ so that the summation runs from $i = 1$ to $n_1$.  
Similar to the expression for $\widetilde{\eta}_1$ in \eqref{36q}, we define
\begin{equation}
\widetilde{\eta}_{11}(t,\theta)
=
\sum_{i=1}^{n_1}
\left\{
\sum_{m=0}^{\widetilde{M}_i}
c_{im} Q_m(\theta)
\right\}
e^{\widetilde{\rho}_i t}.
\tag{48}\label{40q}
\end{equation}

Set
$$
\widetilde{z}_{11} = z_1 - \widetilde{\eta}_{11}.
$$

Then, by the same reasoning,
$$
\mathcal{L}\widetilde{z}_{11}
=
O(e^{\widetilde{\rho}_{n_1+1} t})
=
O(e^{3\rho_1 t}).
$$

Combining this with \eqref{33q} and \eqref{34q}, we obtain
$$
z_{11}=O\left(e^{\tilde{\rho}_1 t}\right)=O\left(e^{2 \rho_1 t}\right)
$$

Since no $\rho_i$ lies between $\widetilde{\rho}_1$ and $\widetilde{\rho}_{r_1+1}$,  
Lemma A.8(ii) of \cite{4z} implies that
$$
z_{11} = O(e^{3\rho_1 t}).
$$

We repeat a procedure similar to Step 2,  
replacing $\widetilde{\rho}_1 = 2\rho_1$ by $\widetilde{\rho}_{r_1+1} = 3\rho_1$.  
If $\rho_{r_1+1} < 4\rho_1$, we imitate the argument of Case 1;  
if $\rho_{r_1+1} > 4\rho_1$, we imitate the argument of Case 2,  
selecting the largest integer $n_2$ such that  
$\widetilde{\rho}_{n_2} < 4\rho_1$.  
Repeating this procedure a finite number of times and up to $\widetilde{\rho}_{l_1}$.

\textbf{Step 3.}  
As in Step 1, we now replace $\widetilde{\rho}_1$ with $\widetilde{\rho}_{l_1+1}$,  
and  $1,\, r_1,\,$ and $1$  
with $r_1+1,\, r_2,\,$ and $l_1+1$.  
Since $\rho_{r_2} < \widetilde{\rho}_{l_1+1}$,  
it follows from \eqref{35q} and Lemma A.8(ii) of \cite{4z}that
$$
z_1(t,\theta)
=
\sum_{i = r_1+1}^{r_2}
c_i Q_i(\theta) e^{\rho_i t}
+
O(e^{\widetilde{\rho}_{l_1+1} t}).
$$

we can discard the terms $e^{\rho_i t}$ for $i = 1,\ldots,r_1$ since $
\widetilde{z}_1 = O(e^{\widetilde{\rho}_1 t}).
$
Define
$$
\eta_2(t,\theta)
=
\sum_{i = r_1+1}^{r_2}
c_i X_i(\theta) e^{\rho_i t},
$$
and set
$$
z_2 = \widetilde{z}_1 - \eta_2.
$$

Then $\mathcal{L}\eta_2 = 0$ and  
$$
z_2 = z - \eta_1 - \widetilde{\eta}_1 - \eta_2,
\qquad
z_2 = O(e^{\widetilde{\rho}_{l_1+1} t}).
$$

\textbf{Step 4.}  
We proceed similarly to Step 2.  
Assume a function $\widetilde{\eta}_2$ is chosen and define
\begin{equation}
\widetilde{z}_2 = z_2 - \widetilde{\eta}_2.
\tag{49}\label{41q}
\end{equation}
Then
$$
\mathcal{L} \tilde z_2=f(z)-\mathcal{L} \tilde{\eta}_1-\mathcal{L} \tilde{\eta}_2
$$
Note that
$$
f(z)
=
f(z_2 + \eta_1 + \widetilde{\eta}_1 + \eta_2)
=
\sum_{i=2}^{\infty} c_i \left( z_2 + \eta_1 + \widetilde{\eta}_1 + \eta_2 \right)^i.
$$

We analyze
$
\sum_{i=2}^{\infty}
c_i \left( \eta_1 + \widetilde{\eta}_1 + \eta_2 \right)^i
$ as in Step 2. Recall from Step 2 and \eqref{36q} that by choosing $\widetilde{\eta}_1$ appropriately,  
we used $\mathcal{L} \widetilde{\eta}_1$ to cancel the terms
$e^{\widetilde{\rho}_i t}$ in $f(z)$ for $i = 1,\ldots,l_1$.  
Proceeding similarly, we can find a function $\widetilde{\eta}_2$ of the form
$$
\widetilde{\eta}_2(t,\theta)
=
\sum_{i = l_1+1}^{l_2}
\left\{
\sum_{m=0}^{\widetilde{M}_i}
c_{im} Q_m(\theta)
\right\}
e^{\widetilde{\rho}_i t}
$$
to cancel the terms $e^{\widetilde{\rho}_i t}$ in $f(z)$ for $i = l_1+1,\ldots,l_2$.

As indicated in \eqref{41q}, defining $\widetilde{z}_2$ accordingly, we obtain
$$
\mathcal{L} z_2
=
O\left( e^{\widetilde{\rho}_{l_2+1} t} \right).
$$

Repeating the argument above completes the proof in the case
$\mathcal{I}_\rho \cap \mathcal{I}_{\widetilde{\rho}} = \emptyset$.

We now consider the more general situation in which some $\rho_i$ can be expressed
as a positive-integer linear combination of $\rho_1,\ldots,\rho_{i-1}$.  
Exponential terms in $t$ will appear in the solution of the equation
$L_i \phi_i = a_i$ when a value $\rho_i$ coincides with a certain $\widetilde{\rho}_{i'}$.   
According to Lemma A.2 of \cite{4z}, such exponential items will appear
one after another during the iteration process.

As a concrete example, we assume $\rho_{r_1} = \widetilde{\rho}_1$,
replacing the strict inequality in \eqref{28q}.
This is the first instance where some $\rho_i$ becomes equal to a
$\widetilde{\rho}_{i'}$.

In this case, we still have
$$
z = O(e^{\widetilde{\rho}_1 t})
$$
and
$$
\mathcal{L}z
=
O(e^{2\rho_1 t})
=
O(e^{\widetilde{\rho}_1 t}).
$$

Following the same procedure as in Step 1,  
choose an index $r_* \in \{1,\ldots,r_1-1\}$ such that
$$
\rho_{r_*}
<
\rho_{r_*+1}
=
\cdots
=
\rho_{r_1}
=
\widetilde{\rho}_1
=
2\rho_1.
$$

By Lemma A.8(ii) of \cite{4z}, we obtain
$$
z(t,\theta)
=
\sum_{i=1}^{r_*}
c_i Q_i(\theta) e^{\rho_i t}
+
O\!\left(t e^{\widetilde{\rho}_1 t}\right).
$$

Define
$$
\eta_1(t,\theta)
=
\sum_{i=1}^{r_*}
c_i X_i(\theta) e^{-\rho_i t}.
$$

With $z_1$ defined similarly to \eqref{32q}, we have
$$
z_1
=
O\!\left(t e^{\widetilde{\rho}_1 t}\right).
$$

Next, proceeding similarly to Step 2, for each
$1 \le i \le l_1$ and $0 \le m \le \widetilde{M}_i$,
we solve equation \eqref{38q}.  
If $\rho_m \ne \widetilde{\rho}_i$, then
$\widetilde{\eta}_{im}(t)$ has the same expression as in \eqref{39q}.  
If $\rho_m = \widetilde{\rho}_i$, then $\widetilde{\eta}_{im}$ takes the form
\begin{equation}
\widetilde{\eta}_{im}(t)
=
c_{i1m}\, t\, e^{\widetilde{\rho}_i t}
+
c_{i0m}\, e^{\widetilde{\rho}_i t}.
\tag{50}\label{42q}
\end{equation}

Using the definition of $\widetilde{\eta}_1$ in \eqref{36q},
the new expression for $\widetilde{\eta}_{im}(t)$ in \eqref{42q},
and $\widetilde{z}_1$ in \eqref{34q},
we obtain \eqref{35q}.  
Repeating the same argument as above yields the desired result.

Now we discard multiple numbers and define a new index sequence
$\{\mu_i\}_{i \ge 1}$.  
Clearly,
$\mu_1 = \rho_1 = 1$ and
$\mu_2 = \min\{ 2\rho_1,\; \rho_{n+1} \}$.

We set
$$
\phi_m(t,\theta)
=
\sum_{\rho_i \le \mu_m}
c_i Q_i(\theta) e^{\rho_i t},
$$
and
$$
\widetilde{\phi}_m(t,\theta)
=
\sum_{\widetilde{\rho}_i \le \mu_m}
\sum_{j=0}^{i}
\left\{
\sum_{l=0}^{\widetilde{M}_i}
c_{ijl} Q_l(\theta)
\right\}
t^j e^{\widetilde{\rho}_i t}.
$$

We note that $\phi_m$ is a solution of $\mathcal{L}\phi_m = 0$ and $\widetilde{\phi}_m$ arises from the nonlinear term $f(z)$.  
In the special case that
$\mathcal{I}_\rho \cap \mathcal{I}_{\widetilde{\rho}} = \emptyset$,
$$
\widetilde{\phi}_m(t,\theta)
=
\sum_{\widetilde{\rho}_i \le \mu_m}
\left\{
\sum_{l=0}^{\widetilde{M}_i}
c_{il} Q_l(\theta)
\right\}
e^{\widetilde{\rho}_i t}.
$$

This completes the proof of (29) in Theorem 3.1.
We now discuss the remaining cases.  
We will not present the full details, and for convenience of notation,
the sequence $\{\rho_i\}_{i\ge 1}$ is again understood without multiple numbers.  
Accordingly, the index sets $\mathcal{I}$, $\mathcal{I}_{\rho}$, 
and $\mathcal{I}_{\widetilde{\rho}}$ can be defined in the same way  
(note that these indices differ from those used earlier).

The ordering of the sequences $\{\rho_i\}_{i\ge 1}$ and 
$\{\widetilde{\rho}_i\}_{i\ge 1}$ becomes
\begin{equation}
\rho_1 < \cdots < \rho_{r_1}
\le \widetilde{\rho}_1
< \cdots <
\widetilde{\rho}_{l_1}
\le \rho_{r_1+1}
< \cdots <
\rho_{r_2}
\le \widetilde{\rho}_{l_1+1}
< \cdots .
\tag{51}\label{52}
\end{equation}

For case \eqref{30}, we set
\[
\begin{aligned}
\widetilde{M}_i
=
\max
\{ &
 k n_1 + (k+1)n_2 + \cdots + (k + i_1 - 1)n_{r_1}
: \\ &
n_1,\ldots,n_{r_1} \text{ are nonnegative integers satisfying \eqref{29q}}
\}.
\end{aligned}
\]

For case \eqref{34}, we set
$$
\begin{aligned}
\widetilde{M}_i
=
\max
\{ &
k n_1 + (k+1)n_2 + \cdots + j n_{j-k+1}
+ 0\, n_{j-k+2}
+ (j+1)n_{j-k+3}
+ \cdots
+ (r_1 + k - 2)n_{r_1}
:\; \\&
n_1,\ldots,n_{r_1} \text{ are nonnegative integers satisfying \eqref{29q}}
\}.
\end{aligned}
$$

For case \eqref{35}, we set
$$
\widetilde{M}_i
=
\max
\left\{
k n_1 + (k+1)n_2 + \cdots + (r_1 + k - 1)n_{r_1}
:\;
n_1,\ldots,n_{r_1} \text{ are nonnegative integers satisfying \eqref{29q}}
\right\}.
$$

For case \eqref{36}, we set
$$
\widetilde{M}_i
=
\max
\left\{
0\, n_1
+ k n_2
+ \cdots
+ (r_1 + k - 2)n_{r_1}
:\;
n_1,\ldots,n_{r_1} \text{ are nonnegative integers satisfying \eqref{29q}}
\right\}.
$$

All remaining steps follow the same pattern as in the proof of \eqref{30}.  
Thus Theorem 3.1 is now proved in full.
\end{proof}
\textbf{Remark.}
When $F \neq 0$, one may rewrite the equation in the form
$\mathcal{L}\!\left(z + C e^{\frac{2p-\alpha}{p+1} t}\right)
= f(z)$
for a suitable constant $C$.  
Consequently, the expansion will contain terms involving
$e^{\frac{2p-\alpha}{p+1} t}$.

\section{Existence of Solutions}
To establish the existence of solutions to \eqref{3}, we introduce a weighted H\"older space.  
In this space, we apply the contraction mapping principle so that the fixed-point leads the existence of a solution to \eqref{3}.  
Define
$$
\|v\|_{C_\mu^i\left((-\infty,t_0] \times S^{N-1}\right)}
=
\sum_{j=0}^i
\sup_{(t,\theta)\in(-\infty,t_0]\times S^{N-1}}
e^{-\mu t}
\left| \nabla^j v(t,\theta) \right|,
$$
and
$$
\|v\|_{C_\mu^{i,a}\left((-\infty,t_0] \times S^{N-1}\right)}
=
\|v\|_{C_\mu^i\left((-\infty,t_0] \times S^{N-1}\right)}
+
\sup_{t \le t_0 - 1}
e^{-\mu t}
\left[
\nabla^i v
\right]_{C^a\!\left([t-1,t+1] \times S^{N-1}\right)},
$$
where $[\cdot]_{C^a}$ means the H\"older seminorm.

\textbf{Definition.}
The weighted H\"older space
$C_\mu^{i,\alpha}\left((-\infty,t_0] \times S^{N-1}\right)$  
consists of those functions $v \in C^i\left((-\infty,t_0] \times S^{N-1}\right)$
for which the norm
$\|v\|_{C_\mu^{i,\alpha}\left((-\infty,t_0] \times S^{N-1}\right)}$
is finite.

Let $\mathcal{L}$ be the linear operator defined in \eqref{eq:L} and let $\mu>0$.  
For a function
$g \in C_\mu^{0,\alpha}\left((-\infty,t_0] \times S^{N-1}\right)$,  
we consider the linear equation
\begin{equation}
\tag{52}\label{eq:eqe}
   \mathcal{L} v = g.
\end{equation}

We shall impose a suitable boundary condition at $t = t_0$ so that
$$
\mathcal{L} :
C_\mu^{2,\alpha}\!\left((-\infty,t_0] \times S^{N-1}\right)
\;\longrightarrow\;
C_\mu^{0,\alpha}\!\left((-\infty,t_0] \times S^{N-1}\right)
$$
admits a bounded inverse.  
We begin with the problem
\begin{equation}
\tag{53}\label{eq:Lv-problem}
\begin{cases}
\mathcal{L} v = g, & \text{in } (-\infty,t_0)\times S^{N-1},\\[2mm]
v = \varphi,       & \text{on } \{t_0\}\times S^{N-1}.
\end{cases}
\end{equation}

\begin{lem}
Let $\mu > 0$,  
$g \in C_\mu^0\!\left((-\infty,t_0] \times S^{N-1}\right)$,
and $\varphi \in C^0(S^{N-1})$.  
Then the above problem admits \emph{at most one} solution  
$v \in C_\mu^2\!\left((-\infty,t_0] \times S^{N-1}\right)$.
\end{lem}

\begin{proof}
We give the argument for the case
$\delta^{(k_0-1)} < \alpha < \delta^{(k_0)}$.
Let $g = 0$ and $\varphi = 0$, suppose
$v \in C_\mu^2\!\left((-\infty,t_0] \times S^{N-1}\right)$
is a solution of the problem.  
For each $k \ge 0$, define
$$
v_k(t)
=
\int_{S^{N-1}} v(t,\theta)\, Q_k(\theta)\, d\theta.
$$

Then $\mathcal{L}_k(v_k) = 0$ and $v_k(t_0) = 0$.  
Hence $v_k$ is a linear combination of elements in
$\operatorname{Ker}(\mathcal{L}_k)$.
For $k < k_0$, we have
\[
v_k(t)
=
c_k^1 e^{\Re(\sigma_1^{(k)}) t} \cos(\gamma t)
+
c_k^2 e^{\Re(\sigma_1^{(k)}) t} \sin(\gamma t),
\]
or
\[
c_k^1 e^{\sigma_1^{(k)} t}
+
c_k^2 t\, e^{\sigma_1^{(k)} t},
\]
or
\[
c_k^1 e^{\sigma_1^{(k)} t}
+
c_k^2 e^{\sigma_2^{(k)} t},
\]
(where $\Re(\sigma_1^{(k)}) < 0$ and $\sigma_1^{(k)}, \sigma_2^{(k)} \le 0$).

From the assumption that $\lim_{t\to -\infty} v_k(t) = 0$ and $v_k(t_0)=0$,  
it follows immediately that $v_k(t) \equiv 0$.

When $k \ge k_0$, we have
$$
v_k(t)
=
c_k^1 e^{\sigma_1^{(k)} t}
+
c_k^2 e^{\sigma_2^{(k)} t}.
$$

Obtain by Green identity that
$$
\int_{-\infty}^{t_0}
\left[
(\partial_t v_k)^2
-
\frac{1}{2}
\left(
N-2 + \frac{2\alpha+4}{p+1}
\right)
(v_k^2)_t
+
\Big(
\lambda_k - (\alpha+2)\big( N-2 + \tfrac{\alpha+2}{p+1} \big)
\Big)
v_k^2
\right]
dt = 0.
$$

Hence,
$$
\int_{-\infty}^{t_0}
\left[
(\partial_t v_k)^2
+
\Big(
\lambda_k - (\alpha+2)\big( N-2 + \tfrac{\alpha+2}{p+1} \big)
\Big)
v_k^2
\right]
dt = 0.
$$

It's clear that the coefficient
$
\lambda_k - (\alpha+2)\left( N-2 + \frac{\alpha+2}{p+1} \right)
> 0
$
when $k \ge k_0$,
thus $v_k(t) \equiv 0$,  
$v \equiv 0$.
\end{proof}
We now estimate the $C^{2,\alpha}$ norm of the solution to \eqref{eq:Lv-problem}.

\begin{lem}
Let $\alpha \in (0,1)$, $\mu > 0$,
$g \in C_\mu^{0,\alpha}\!\left((-\infty,t_0] \times S^{N-1}\right)$
and $\varphi \in C^{2,\alpha}(S^{N-1})$.
Suppose 
$v \in C_\mu^{2,\alpha}\!\left((-\infty,t_0] \times S^{N-1}\right)$
is a solution of \eqref{eq:Lv-problem}.  
Then
\begin{equation}
\tag{54}\label{qs}
\|v\|_{C_\mu^{2,\alpha}\left((-\infty,t_0] \times S^{N-1}\right)}
\le
C\Big[
\|v\|_{C_\mu^0\left((-\infty,t_0] \times S^{N-1}\right)}
+
\|g\|_{C_\mu^{0,\alpha}\left((-\infty,t_0] \times S^{N-1}\right)}
+
e^{-\mu t_0}\|\varphi\|_{C^{2,\alpha}(S^{N-1})}
\Big],
\end{equation}
where $C>0$ depends only on $N,\alpha,\mu$ and is independent of $t_0$.
\end{lem}

\begin{proof}
 Fix a $t \le t_0$ and consider two cases, (i) $t < t_0 - 2$, 
(ii) $t_0 - 2 \le t \le t_0$. They are proved by the interior and boundary Schauder estimate respectively. 

For case (i), By the interior Schauder estimate,
\[
\begin{aligned}
&\sum_{j=0}^2 \sup_{S^{N-1}} |\nabla^j v(t,\cdot)|
+
\left[\nabla^2 v\right]_{C^\alpha([t-1,t+1]\times S^{N-1})}
\\
&\qquad\le
C\Big[
\|v\|_{L^\infty([t-2,t+2]\times S^{N-1})}
+
\|g\|_{L^\infty([t-2,t+2]\times S^{N-1})}
+
[g]_{C^\alpha([t-2,t+2]\times S^{N-1})}
\Big],
\end{aligned}
\]
where $C>0$ is independent of $t$.

To estimate the H\"older seminorm of $g$ on $[t-2,t+2]\times S^{N-1}$,
take arbitrary
$(t_1,\theta_1),(t_2,\theta_2)\in[t-2,t+2]\times S^{N-1}$  
with $(t_1,\theta_1)\neq(t_2,\theta_2)$.
Then we split into two subcases:
$|t_1 - t_2| \le 2$ and $|t_1 - t_2| > 2$.
We present the proof of the first subcase and the other follows directly. 

For $|t_1 - t_2| \le 2$, there exists some
$t' \in [t-1,t+1]$ such that  
$t_1,t_2 \in [t'-1,t'+1]$.  
Hence,
\[
[g]_{C^\alpha([t-2,t+2]\times S^{N-1})}
\le
\max\Big\{
\sup_{t' \in [t-1,t+1]}
[g]_{C^\alpha([t'-1,t'+1]\times S^{N-1})},
\ 
\|g\|_{L^\infty([t-2,t+2]\times S^{N-1})}
\Big\}.
\]
Thus
$$
\begin{aligned} & \sum_{j=0}^2 \sup _{S^{N-1}}\left|\nabla^j v(t, \cdot)\right|+\left[\nabla^2 v\right]_{C^\alpha\left([t-1, t+1] \times S^{N-1}\right)} \\ & \quad \leq C\left[\|v\|_{L^{\infty}\left([t-2, t+2] \times S^{N-1}\right)}+\|g\|_{L^{\infty}\left([t-2, t+2] \times S^{N-1}\right)}+\sup _{t^{\prime} \in[t-1, t+1]}[g]_{C^\alpha\left(\left[t^{\prime}-1, t^{\prime}+1\right] \times S^{N-1}\right)}\right].\end{aligned}
$$

Multiplying both sides by $e^{-\mu t}$ and then taking the supremum over
$t \in (-\infty, t_0 - 2)$, we obtain
$$
\begin{aligned}
& \sum_{j=0}^2
\sup_{t \in (-\infty, t_0 - 2)}
\sup_{S^{N-1}}
e^{-\mu t} \left| \nabla^j v(t,\cdot) \right|
+
\sup_{t \in (-\infty, t_0 - 2)}
e^{-\mu t}
\left[
\nabla^2 v
\right]_{C^\alpha([t-1,t+1] \times S^{N-1})}
\\
&\qquad \le
C\left[
\|v\|_{C_\mu^0((-\infty,t_0] \times S^{N-1})}
+
\|g\|_{C_\mu^{0,\alpha}((-\infty,t_0] \times S^{N-1})}
\right],
\end{aligned}
$$
where the constant $C>0$ is independent of $t_0$. 

For case (ii), it's similar with case (i) as we know by the boundary Schauder estimate that 
\begin{equation*}
\begin{aligned} & \sum_{j=0}^2 \sup _{S^{N-1}}\left|\nabla^j v(t, \cdot)\right|+\left[\nabla^2 v\right]_{C^\alpha\left(\left[t_0-3, t_0\right] \times S^{N-1}\right)} \\ & \quad \leq C\left[\|v\|_{L^{\infty}\left(\left[t_0-4, t_0\right] \times S^{N-1}\right)}+\|g\|_{L^{\infty}\left(\left[t_0-4, t_0\right] \times S^{N-1}\right)}+[g]_{C^\alpha\left(\left[t_0-4, t_0\right] \times S^{N-1}\right)}+\|\varphi\|_{C^{2, \alpha}\left(S^{N-1}\right)}\right],\end{aligned}
\end{equation*}
we omit the details here. Thus we get the desired estimate by case (i) and case (ii). 
\end{proof}

\begin{lem}
Assume $\mu > 0$ and  $\mu \neq \sigma_1^{(k)}$ for every $k \ge 1$.
Let $T$ and $t_0$ be constants with $t_0 \le 0$ and $T - t_0 \le -4$, suppose
$g \in C^0([T,t_0] \times S^{N-1})$.
$v \in C^2([T,t_0] \times S^{N-1})$ satisfies
$$
\begin{cases}
\mathcal{L} v = g & \text{in } (T,t_0)\times S^{N-1},\\
v = 0             & \text{on } (\{T\} \cup \{t_0\})\times S^{N-1},
\end{cases}
$$
and for every $t \in [T,t_0]$ and every $k = 0,1,2,\ldots,K,$ 
$$
\int_{S^{N-1}} v(t,\theta)\, Q_k(\theta)\, d\theta = 0,
$$
where $K$ is the largest integer such that $\sigma_1^{(K)} < \mu$,  
then
$$
\sup_{(t,\theta)\in [T,t_0]\times S^{N-1}}
e^{-\mu t} |v(t,\theta)|
\le
C
\sup_{(t,\theta)\in [T,t_0]\times S^{N-1}}
e^{-\mu t} |g(t,\theta)|,
$$
where $C>0$ depends only on $N$ and $\mu$ and independent of $T$ and $t_0$.
\end{lem}

\begin{proof}
Assume the conclusion is false.
Then there exist sequences
$\{T_i\}$, $\{t_i\}$, $\{v_i\}$ and $\{g_i\}$  
with $t_i \le 0$ and $T_i - t_i \le -4$, such that
$$
\begin{cases}
\mathcal{L} v_i = g_i & \text{in } (T_i,t_i)\times S^{N-1},\\
v_i = 0               & \text{on } (\{T_i\} \cup \{t_i\})\times S^{N-1},
\end{cases}
$$
and 
$$
\sup_{(t,\theta)\in [T_i,t_i]\times S^{N-1}}
e^{-\mu t} |g_i(t,\theta)| = 1,
\qquad
\sup_{(t,\theta)\in [T_i,t_i]\times S^{N-1}}
e^{-\mu t} |v_i(t,\theta)|
\to \infty.\quad (\text{as } i \to \infty.)
$$

From each interval $(T_i,t_i)$ ,choose a point $t_i^* \in (T_i,t_i)$ such that
$$
M_i
=
\sup_{S^{N-1}}
e^{-\mu t_i^*}
\left| v_i(t_i^*,\cdot) \right|
=
\sup_{(t,\theta)\in [T_i,t_i]\times S^{N-1}}
e^{-\mu t} |v_i(t,\theta)|.
$$

Then $M_i \to \infty$ as $i \to \infty$.
Define 
$$
\widetilde{v}_i(t,\theta)
=
M_i^{-1} e^{-\mu t_i^*}
v_i(t + t_i^*, \theta),
$$
and
$$
\widetilde{g}_i(t,\theta)
=
M_i^{-1} e^{-\mu t_i^*}
g_i(t + t_i^*, \theta).
$$

Then
$$
\sup_{S^{N-1}}
|\widetilde{v}_i(0,\cdot)|
= 1,
$$

For every $(t,\theta) \in [T_i - t_i^*,\, t_i - t_i^*] \times S^{N-1}$ we have
\begin{equation}
\tag{55}\label{eq:wqw}
\left|e^{-\mu t} \tilde{v}_i(t, \theta)\right| \leq 1 .
\end{equation}
   
Furthermore, for all
$t \in [T_i - t_i^*,\, t_i - t_i^*] \times S^{N-1}$, we have
$$
\mathcal{L}\widetilde{v}_i = \widetilde{g}_i.
$$

Passing to a subsequence, we may assume that there exist
$\tau_- \in \mathbb{R}^- \cup \{-\infty\}$ and
$\tau_+ \in \mathbb{R}^+ \cup \{\infty\}$
such that
$$
T_i - t_i^* \to \tau_-,
\qquad
t_i - t_i^* \to \tau_+.
$$

In fact, we obtain from \eqref{eq:wqw} that
$$
|\widetilde{v}_i|
\le
C e^{\mu (t_i^* - T_i)}
\qquad\text{on } (T_i - t_i^*,\, T_i - t_i^* + 2)\times S^{N-1},
$$
hence
$$
\left|
\frac{d^2 \widetilde{v}_i}{dt^2}
+
\left( N-2 + \frac{4}{p+1} \right)
\frac{d\widetilde{v}_i}{dt}
+
\Delta_\theta \widetilde{v}_i
\right|
\le
C e^{\mu (t_i^* - T_i)}
\qquad\text{on } (T_i - t_i^*,\, T_i - t_i^* + 2)\times S^{N-1}.
$$

Since $\widetilde{v}_i = 0$ on $\{T_i - t_i^*\} \times S^{N-1}$, it follows that
$$
|\nabla \widetilde{v}_i|
\le
C e^{\mu (t_i^* - T_i)}
\qquad\text{on } (T_i - t_i^*,\, T_i - t_i^* + 1)\times S^{N-1},
$$
so
\begin{equation*}
\begin{aligned}\left|\tilde{v}_i(t, \theta)\right| & =\left|\tilde{v}_i(0, \theta)-\tilde{v}_i\left(T_i-t_i^*, \theta\right)\right| \\ & \leqslant \sup _{t \in\left(T_i-t_i^*, T_i-t_i^*+1\right)}\left|\partial_t \tilde{v}_i(t, \theta)\right|\left|T_i-t_i^*\right| \\ & \leqslant \sup _{(t, \theta) \in\left(T_i-t_i^*, T_i-t_i^*+1\right) \times S^{N-1}}\left|\nabla \tilde{V}_i(t, \theta)\right| \cdot\left|T_i-t_i^*\right| \\ & \leqslant c e^{\mu\left(t_i^*-T_i\right)}\left|T_i-t_i^*\right|\end{aligned}
\end{equation*}
Consequently $T_i - t_i^*$ keeps a definite distance away from $0$.
Similarly, $t_i - t_i^*$ also keeps a definite distance away from $0$.
Hence,
$$
0 \in (\tau_-, \tau_+).
$$

Assume
$$
\widetilde{v}_i \to \widehat{v}
\quad\text{uniformly on compact subsets of } (\tau_-, \tau_+) \times S^{N-1}.
$$

Also $\widetilde{g}_i$ converges uniformly to zero
on every compact subset of $(\tau_-, \tau_+)\times S^{N-1}$.
Therefore, $\widehat{v} \neq 0$, and
\begin{equation}
\tag{56}\label{eq:wqwq}
|e^{-\mu t} \widehat{v}(t,\theta)| \le 1 \quad
\text{for all}
(t,\theta) \in (\tau_-, \tau_+) \times S^{N-1}
\end{equation}
 Also we have
$$
\mathcal{L}\widehat{v} = 0 \quad on ~ (\tau_-, \tau_+) \times S^{N-1},
$$
and
\begin{equation}
\tag{57}\label{eq:2121}
\lim_{t \to \tau} \widehat{v}(t,\theta) = 0.
\end{equation}

Where $\tau=\tau_{-}$ or $\tau_{+}$ if it is finite.

Now, for any $k \ge 0$, define
$$
\widehat{v}_k(t)
=
\int_{S^{N-1}}
\widehat{v}(t,\theta)\, Q_k(\theta)\, d\theta.
$$

Then $\mathcal{L}_k(\widehat{v}_k)=0$, $\widehat{v}_k$ is a linear
combination of a basis of $\operatorname{Ker}(\mathcal{L}_k)$.
Choose $k$ such that $\sigma_1^{(k)} > \mu > 0$.  
Then
$$
\widehat{v}_k(t)
=
c_k^1 e^{\sigma_1^{(k)} t}
+
c_k^2 e^{\sigma_2^{(k)} t}.
$$

From \eqref{eq:wqwq} we know 
$$
| e^{-\mu t} \widehat{v}_k(t) | \le C \quad for ~ all ~t \in (\tau_-,\tau_+).
$$

If $\tau_+ = +\infty$, then  $c_k^1 = 0$,  thus
$\widehat{v}_k(t) = c_k^2 e^{\sigma_2^{(k)} t}$
which decays exponentially as $t \to +\infty$.
If $\tau_+$ is finite, then
$\lim_{t\to \tau_+} \widehat{v}_k(t) = 0$ by \eqref{eq:2121}.

Similarly, if $\tau_- = -\infty$, then $c_k^2 = 0$,  hence
$\widehat{v}_k(t) = c_k^1 e^{\sigma_1^{(k)} t}$
which decays exponentially as $t \to -\infty$.
If $\tau_-$ is finite, then 
$\lim_{t\to \tau_-} \widehat{v}_k(t) = 0$  by \eqref{eq:2121}.

Therefore,
$$
\int_{\tau_-}^{\tau_+}
\left[
(\partial_t \widehat{v}_k)^2
+
\left(
\lambda_k
-
(\alpha+2)\left(
N-2 + \frac{\alpha+2}{p+1}
\right)
\right)
\widehat{v}_k^2
\right] dt = 0.
$$

We have
$$
\lambda_k
-
(\alpha+2)\left(
N-2 + \frac{\alpha+2}{p+1}
\right)
> 0.
$$
Since $\sigma_1^{(k)} > \mu > 0$.
So $\widehat{v}_k \equiv 0$
for every $k$ satisfying $\sigma_1^{(k)} > \mu$ .

From the assumption that
$$
\int_{S^{N-1}} \widetilde{v}_i(t,\theta) Q_k(\theta)\, d\theta = 0
\qquad
\text{for all } k = 0,1,\ldots,K
\text{ with }
\sigma_1^{(K)} < \mu,
$$
 since $\widetilde{v}_i \to \widehat{v}$ as $i\to \infty$, we deduce that
$$
\widehat{v}_k \equiv 0
\qquad
\text{for all } k = 0,1,\ldots,K \text{ with }
\sigma_1^{(K)} < \mu.
$$

Combining the above, we conclude that $\widehat{v}_k \equiv 0$ for all $k \ge 0$,
therefore $\widehat{v} \equiv 0$ and it contradicts with our earlier conclusion that
$\widehat{v} \neq 0$. Lemma 4.3 is proved.
\end{proof}
\begin{lem}
Let $\alpha \in (0,1)$ and $\mu > \sigma_1^{(K)}$ (with $K \ge 0$), and let
$g \in C_\mu^{0,\alpha}\!\left((-\infty,t_0] \times S^{N-1}\right)$
satisfying that
$g(t,\cdot) \in \operatorname{span}\{Q_0, Q_1, \ldots, Q_K\}$
for all $t \le t_0$.
Then equation \eqref{eq:eqe} admits a unique solution
$v \in C_\mu^{2,\alpha}\!\left((-\infty,t_0] \times S^{N-1}\right)$,
and for every $t \le t_0$,
$
v(t,\cdot) \in \operatorname{span}\{Q_0, Q_1, \ldots, Q_K\}.
$
Moreover, the map $g \mapsto v$ is linear, and
$$
\| v \|_{C_\mu^{2,\alpha}\left((-\infty,t_0] \times S^{N-1}\right)}
\le
C \| g \|_{C_\mu^{0,\alpha}\left((-\infty,t_0] \times S^{N-1}\right)},
$$
where the constant $C>0$ depends only on $N,\alpha,\mu$ and is independent of $t_0$.
\end{lem}
\begin{proof}
For each $k = 0,1,\ldots,K$, define
$$
g_k(t)
=
\int_{S^{N-1}}
g(t,\theta)\, Q_k(\theta)\, d\theta.
$$

Then
$$
\| g_k \|_{C_\mu^{0,\alpha}((-\infty,t_0])}
\le
C \| g \|_{C_\mu^{0,\alpha}((-\infty,t_0] \times S^{N-1})},
$$
and
$$
g(t,\theta)
=
\sum_{k=0}^K g_k(t)\, Q_k(\theta).
$$

Let $\mathcal{L}_k$ denote the linear operator from \eqref{wqwqq}.
Consider the ordinary differential equation
\begin{equation}
\tag{58}\label{qw}
\mathcal{L}_k v_k = g_k.
\end{equation}
Then
\begin{equation}
v_k(t)=B_k^1 \int_{-\infty}^t e^{\sigma_1^{(k)}(t-s)} g_k(s) d s-B_k^2 \int_{-\infty}^t e^{\sigma_2^{(k)}(t-s)} g_k(s) d s,
\tag{59}\label{61}
\end{equation}
or
\begin{equation*}
\begin{aligned}
v_k(t) &= {B_k^1}^{*}  \int_{-\infty}^t e^{R(\sigma_1^{(k)})(t-s)} g_k(s) \sin\bigl(I(\sigma_1^{(k)}) (t-s)\bigr) \, ds ,\\
\end{aligned}
\end{equation*}
or
\begin{equation*}
\begin{aligned}
v_k(t)=e^{\sigma_1^{(k)} t} \int_{-\infty}^t \int_{-\infty}^s e^{-\sigma_1^{(k)} \gamma} g_k(\gamma) d \gamma d s,
\end{aligned}
\end{equation*}
where $B_k^1~and~{B_k^1}^{*} $ are constants.
We now write down the argument for the first and the other are similar.

A direct computation shows that, for all $t \le t_0$,
$$
\begin{gathered}
e^{-\mu t} |v_k(t)|
\le
C \sup_{t \le t_0} e^{-\mu t} |g_k(t)|
=
C \| g_k \|_{C_\mu^0((-\infty,t_0])},
\\[2mm]
e^{-\mu t}\left( |v_k'(t)| + |v_k''(t)| \right)
\le
C \| g_k \|_{C_\mu^0((-\infty,t_0])}.
\end{gathered}
$$

To estimate the H\"older seminorm of $v_k''$, decompose
$$
v_k''(t) = R_1(t) + R_2(t),
$$
where
$$
R_1(t)
=
B_k^1 (e^{\sigma_1^{(k)} t})'' \int_{-\infty}^t e^{-\sigma_1^{(k)} s} g_k(s)\, ds
-
B_k^1 (e^{\sigma_2^{(k)} t})'' \int_{-\infty}^t e^{-\sigma_2^{(k)} s} g_k(s)\, ds,
$$
and
$$
R_2(t)
=
B_k^1 (e^{\sigma_1^{(k)} t})' e^{-\sigma_1^{(k)} t} g_k(t)
-
B_k^1 (e^{\sigma_2^{(k)} t})' e^{-\sigma_2^{(k)} t} g_k(t).
$$

Thus,
$$
\begin{aligned}
R_1'(t)
= {} &
B_k^1 (e^{\sigma_1^{(k)} t})''' \int_{-\infty}^t e^{-\sigma_1^{(k)} s} g_k(s)\, ds
-
B_k^1 (e^{\sigma_2^{(k)} t})''' \int_{-\infty}^t e^{-\sigma_2^{(k)} s} g_k(s)\, ds
\\
& +
B_k^1 (e^{\sigma_1^{(k)} t})'' e^{-\sigma_1^{(k)} t} g_k(t)
-
B_k^1 (e^{\sigma_2^{(k)} t})'' e^{-\sigma_2^{(k)} t} g_k(t).
\end{aligned}
$$

Similarly, for all $t \le t_0$,
$$
e^{-\mu t} |R_1'(t)|
\le
C \| g_k \|_{C_\mu^0((-\infty,t_0])},
$$
 hence, for all $t \le t_0 - 1$,
$$
\begin{aligned}
e^{-\mu t} [R_1]_{C^\alpha([t-1,t+1])}
&\le
C \| g_k \|_{C_\mu^0((-\infty,t_0])},
\\
e^{-\mu t} [R_2]_{C^\alpha([t-1,t+1])}
&\le
C \| g_k \|_{C_\mu^{0,\alpha}((-\infty,t_0])}.
\end{aligned}
$$

Therefore, for all $t \le t_0 - 1$,
$$
e^{-\mu t}
\left[
v_k''
\right]_{C^\alpha([t-1,t+1])}
\le
C \| g_k \|_{C_\mu^{0,\alpha}((-\infty,t_0])}.
$$

Combining the above estimates obtains that
  \eqref{qw} admits a solution
$v_k \in C_\mu^{2,\alpha}((-\infty,t_0])$
satisfying
\begin{equation}
\tag{60}\label{21212}
||v_k||_{C_\mu^{2,\alpha}((-\infty,t_0])}
\le
C\, ||g_k||_{C_\mu^{0,\alpha}((-\infty,t_0])},
\end{equation}
where the constant $C$ depends only on $N,\alpha,\mu$ and not on $t_0$.

After get the solution of $v_k$ from \eqref{qw} for each $k = 0,1,\ldots,K$, we set
$$
v(t,\theta) = \sum_{k=0}^K v_k(t) Q_k(\theta).
$$

Clearly $\mathcal{L}v = g$, then we can get from \eqref{21212} that
\begin{equation*}
\begin{aligned}
\|v\|_{C_\mu^{2,\alpha}((-\infty,t_0]\times S^{N-1})}
\le
C \sum_{k=0}^K \|v_k\|_{C_\mu^{2,\alpha}((-\infty,t_0])}
\le
C \sum_{k=0}^K \|g_k\|_{C_\mu^{0,\alpha}((-\infty,t_0])}
\le
C \|g\|_{C_\mu^{0,\alpha}((-\infty,t_0]\times S^{N-1})}.
\end{aligned}
\end{equation*}

Thus $v$ is the desired solution.  
Moreover, it's clear that such a solution is unique under the additional condition that  
$
v(t,\cdot) \in \operatorname{span}\{Q_0, Q_1, \ldots, Q_K\}~\text{for every } t \le t_0.
$
\end{proof}
\begin{lem}
Let $\alpha \in (0,1)$, $\mu > 0$ with $\mu \neq \sigma_1^{(k)}$ for every $k \ge 1$, and let
$g \in C_\mu^{0,\alpha}((-\infty,t_0]\times S^{N-1})$ satisfying
$$
\int_{S^{N-1}} g(t,\theta)\, Q_k(\theta)\, d\theta = 0,
\qquad
\forall\, t \le t_0,\quad k = 0,1,\ldots,K,
$$
where $K$ is the largest integer such that $\sigma_1^{(K)} < \mu$.
Then equation \eqref{eq:eqe} admits a unique solution  
$v \in C_\mu^{2,\alpha}((-\infty,t_0]\times S^{N-1})$ satisfying  
$v = 0$ on $\{t_0\} \times S^{N-1}$.
Furthermore,
$$
\|v\|_{C_\mu^{2,\alpha}((-\infty,t_0]\times S^{N-1})}
\le
C\|g\|_{C_\mu^{0,\alpha}((-\infty,t_0]\times S^{N-1})},
$$
where $C>0$ depends only on $N,\alpha,\mu$ and is independent of $t_0$.
\end{lem}

\begin{proof}
Fix any $T \le t_0 - 4$.
We first show that there exists
$v'\in C^{2,\alpha}([T,t_0]\times S^{N-1})$
such that
\begin{equation}
\tag{61}\label{233}
\begin{cases}
\mathcal{L} v' = g & \text{in } (T,t_0)\times S^{N-1},\\[1mm]
v' = 0            & \text{on } (\{T\} \cup \{t_0\})\times S^{N-1}.
\end{cases}
\end{equation}

Observe that the problem \eqref{233} is equivalent to
$$
\begin{cases}\frac{\partial}{\partial t}\left(e^{{A} t} \frac{\partial v'}{\partial t}\right)+e^{{A} t} \Delta_\theta v'+(\alpha+2)\left(N-2+\frac{\alpha+2}{p+1}\right) e^{{A} t} v'=e^{{A} t} g & \text { in }\left(T, t_0\right) \times S^{N-1}, \\ v'=0 & \text { on }\left(\{T\} \cup \{t_0\}\right) \times S^{N-1},\end{cases}
$$
Here ${A} = N - 2 + \frac{2\alpha + 4}{p+1}$.  
Consider the energy functional
$$
\mathcal{G}_T(v')
=
\int_T^{t_0} \int_{S^{N-1}}
\left[
e^{{A} t} (\partial_t v')^2
+
e^{{A} t} |\nabla_\theta v'|^2
-
(\alpha+2)\left( N-2 + \frac{\alpha+2}{p+1} \right) e^{{A} t} {v'}^2
+
2 e^{{A} t} g v'
\right] dt\, d\theta .
$$

Define
$$
\Gamma
=
\left\{
\phi \in H^1(S^{N-1})
:
\int_{S^{N-1}} \phi(\theta)\, Q_k(\theta)\, d\theta = 0,
\quad
k = 0,1,2,\ldots,K,\ 
\sigma_1^{(k)} < \mu
\right\},
$$
set $ \phi=\sum_{k \geq 0} a_k Q_k $, then

$$
\int \phi^2=\sum_{k \geq 0}\left|a_k\right|^2, \quad \int\left|\nabla_\theta \phi\right|^2=\int \phi(-\Delta) \phi=\sum_{k \geq 0} \lambda_k\left|a_k\right|^2 ,
$$
by the definition of $\Gamma$, we know 
$$
a_0=a_1=\cdots=a_K=0,
$$
so
$$
\int\left|\nabla_\theta \phi\right|^2=\sum_{k \geq K+1} \lambda_k\left|a_k\right|^2 \geq \lambda_{K+1} \sum_{k \geq K+1}\left|a_k\right|^2=\lambda_{K+1} \int \phi^2 .
$$

Since $\sigma_{1}^{(K+1)}>\mu >0$, we have
$$
\lambda_{K+1}-
(\alpha+2)\left( N-2 + \frac{\alpha+2}{p+1} \right)>0
,
$$
it follows that for any  
$v' \in H_0^1((T,t_0)\times S^{N-1})$ satisfying  
$v'(t,\cdot) \in \Gamma$ ($t \in (T,t_0)$), we have
$$
\mathcal{G}_T(v')
\ge
\int_T^{t_0} \int_{S^{N-1}}
e^{{A} t} (\partial_t v')^2
+
\left(
\lambda_{K+1}
-
(\alpha+2)\left( N-2 + \frac{\alpha+2}{p+1} \right)
\right)
e^{{A} v'} v^2
+
2 e^{{A} v'} g v'
 dt d\theta .
$$

Hence the functional $\mathcal{G}_T$ is coercive and weakly lower semicontinuous.  
Therefore we may  find a minimizer $v'$ of $\mathcal{G}_T$ in the space
$$
\left\{
v' \in H_0^1((T,t_0)\times S^{N-1})
:
v'(t,\cdot) \in \Gamma \text{ for every } t \in (T,t_0)
\right\}.
$$

Since $g(t,\cdot) \in \Gamma$ for all $t \in (T,t_0)$,  
it follows that $v'$ is a solution of \eqref{233},  
and $v'(t,\cdot) \in \Gamma$ for all $t \in (T,t_0)$.

By Lemma 4.3, we have
$$
\sup_{(t,\theta)\in [T,t_0]\times S^{N-1}}
e^{-\mu t}\, |v'(t,\theta)|
\le
C \sup_{(t,\theta)\in [T,t_0]\times S^{N-1}}
e^{-\mu t}\, |g(t,\theta)|,
$$
where the constant $C>0$ depends only on $N$ and $\mu$ and is independent of $T$ and $t_0$.

For any fixed $T_0 < t_0$, consider the region
$\,[t_0 + T_0,\, t_0] \times S^{N-1} 
\subset  [t_0 + T_0 - 1,\, t_0] \times S^{N-1}$.
By the interior and boundary Schauder estimates, together with the fact that
$v'(t_0,\theta)=0$, we may extract a subsequence 
$v'$ converges 
to a $C^{2,\alpha}$ solution $v$ of \eqref{eq:eqe} on $[t_0 + T_0,\, t_0] \times S^{N-1}$
by Arzel\`a--Ascoli theorem,
with $v=0$ on $\{t_0\}\times S^{N-1}$ as $T\to -\infty$.
By a diagonalization process, we conclude that $v'$ converges to a $C^{2,\alpha}$ solution $v$ of \eqref{eq:eqe} on 
$(-\infty,t_0] \times S^{N-1}$ ,
with $v=0$ on $\{t_0\}\times S^{N-1}$.

Furthermore,
$$
\sup_{(t,\theta)\in [T,t_0]\times S^{N-1}}
e^{-\mu t} |v(t,\theta)|
\le
C \sup_{(t,\theta)\in [T,t_0]\times S^{N-1}}
e^{-\mu t} |g(t,\theta)|,
$$
or
\begin{equation}
\tag{62}\label{qd}
\|v\|_{C_\mu^0((-\infty,t_0]\times S^{N-1})}
\le
C \|g\|_{C_\mu^0((-\infty,t_0]\times S^{N-1})},
\end{equation}
where $C>0$ depends only on $N$ and $\mu$, and not on $t_0$.

Substituting \eqref{qd} into \eqref{qs} (with $\varphi = 0$) and then completes the proof.
\end{proof}
\begin{thm}
Let $\alpha \in (0,1)$ and $\mu > 0$ with $\mu \ne \sigma_1^{(k)}$ for all $k \ge 1$,  let
$g \in C_\mu^{0,\alpha}((-\infty,t_0] \times S^{N-1})$.
Then equation \eqref{eq:eqe} admits a solution
$v \in C_\mu^{2,\alpha}((-\infty,t_0] \times S^{N-1})$ satisfying
$$
\|v\|_{C_\mu^{2,\alpha}((-\infty,t_0] \times S^{N-1})}
\le
C
\|g\|_{C_\mu^{0,\alpha}((-\infty,t_0] \times S^{N-1})},
$$
where the constant $C>0$ depends only on $N,\alpha,\mu$ and is independent of $t_0$.
Moreover, the corresponding map $g \mapsto v$ is linear.
\end{thm}

\begin{proof}
Let $K \ge 0$ be the largest integer such that $\sigma_1^{(K)} < \mu$.  
For $k = 0,1,\ldots,K$, define
$$
g_k(t)
=
\int_{S^{N-1}}
g(t,\theta)\, Q_k(\theta)\, d\theta.
$$

Let
$v_1 \in C_\mu^{2,\alpha}((-\infty,t_0] \times S^{N-1})$
be the unique solution (as in Lemma 4.4) of
$$
\mathcal{L}(v_1)
=
\sum_{k=0}^K g_k(t)\, Q_k(\theta)
\qquad\text{in } (-\infty,t_0] \times S^{N-1}.
$$

Let
$v_2 \in C_\mu^{2,\alpha}((-\infty,t_0] \times S^{N-1})$
to be the unique solution (as in Lemma 4.6) of
$$
\begin{cases}
\mathcal{L} v
=
g - \sum_{k=0}^K g_k Q_k
&\text{in } (-\infty,t_0] \times S^{N-1},
\\[2mm]
v = 0
&\text{on } \{t_0\} \times S^{N-1}.
\end{cases}
$$

Then $v = v_1 + v_2$ is the desired solution.

\end{proof}
\begin{remark}
We denote $\mathcal{L}^{-1}$ the correspondence $g \mapsto v$ in Theorem 4.6.
Then
$$
\mathcal{L}^{-1} :
C_\mu^{0,\alpha}((-\infty,t_0] \times S^{N-1})
\longrightarrow
C_\mu^{2,\alpha}((-\infty,t_0] \times S^{N-1})
$$
is a bounded linear operator, and its operator norm is independent of $t_0$.
\end{remark}
Next we prove the existence of solutions to \eqref{3}.
Let
\[
H(z)
=
z_{tt}
+
\left( N - 2 + 2q \right) z_t
+
(\alpha+2)(N-2+q)\, z
+
\Delta_\theta z
-
f(z)
-
F e^{\frac{2p-\alpha}{p+1} t}.
\]
It's equivalent to prove   $H(z) = 0$ admits a solution.

\begin{thm}
(1) \emph{Radial case.}
Assume $\mu>0$ and $F\neq0$.
Suppose $\hat{z} \in C^{2,\alpha}((-\infty,0])$ satisfying
$$
|\hat{z}(t)| + |\hat{z}'(t)| \to 0
\quad \text{as } t \to -\infty,
$$
and there exists a positive constant $C$ such that
$$
\big|H(\hat{z})\big|
+
\left|\frac{d}{dt} H(\hat{z})\right|
\le
C e^{\mu t}
\quad \forall\, t \in (-\infty,0].
$$
Then there exists $t_0 < 0$ and a solution
$z(t) \in C^{2,\alpha}((-\infty,t_0])$ of ${H}(z) = 0$ such that
$$
|z(t) - \hat{z}(t)|
\le
C e^{\mu t}
\quad \text{for } t \in (-\infty,t_0),
$$
where $C$ is a positive constant.

(2) \emph{Nonradial case.}
Assume $\mu > 0$ and $\mu \ne \sigma_1^{(i)}$ for every $i \ge 1$.
Suppose
$\hat{z} \in C^{2,\alpha}((-\infty,0]\times \mathbb{S}^{N-1})$ satisfies
$$
|\hat{z}(t,\theta)| + |\nabla \hat{z}(t,\theta)| \to 0
\quad \text{as } t \to -\infty
\text{ uniformly for } \theta \in \mathbb{S}^{N-1},
$$
and there exists a positive constant $C$ such that for all
$(t,\theta) \in (-\infty,0]\times \mathbb{S}^{N-1}$,
\begin{equation}
|H(\hat{z})|
+
|\nabla H(\hat{z})|
\le
C e^{\mu t}.
\tag{63}\label{64q}
\end{equation}
Then there exist $t_0 < 0$ and a solution
\[
z(t,\theta)
\in
C^{2,\alpha}\big((-\infty,t_0]\times \mathbb{S}^{N-1}\big)
\quad \text{of } {H}(z) = 0,
\]
such that
$$
|z(t,\theta) - \hat{z}(t,\theta)|
\le
C e^{\mu t}
\quad \text{for } (t,\theta) \in (-\infty,t_0)\times \mathbb{S}^{N-1},
$$
where $C$ is a positive constant.
\end{thm}
\begin{proof}
We only prove case (2) since the radial case (1) is similar and easier.
For any
$\phi \in C_\mu^{2,\alpha}((-\infty,t_0]\times S^{N-1})$, we have
$$
H(\widehat{z} + \phi)
=
H(\widehat{z})
+
\mathcal{L}\phi
-
P(\phi),
$$
where
$$
P(\phi)
=
(\widehat{z} + \Lambda + \phi)^{-p}
-
(\widehat{z} + \Lambda)^{-p}
+
p \Lambda^{-(p+1)} \phi.
$$

Thus $H(\widehat{z} + \phi) = 0$ is equivalent to
\begin{equation}
\mathcal{L}\phi
=
- \big[ H(\widehat{z}) - P(\phi) \big].
\tag{64}\label{hh}
\end{equation}

By Theorem 4.6 and corresponding remark, we may rewrite \eqref{64q} in the form
$$
\phi
=
\mathcal{L}^{-1}\big( -H(\widehat{z}) + P(\phi) \big).
$$

Define the operator
$$
\mathcal{T}(\phi)
=
\mathcal{L}^{-1}\big( -H(\widehat{z}) + P(\phi) \big).
$$
We shall show that, for $t_0 < 0$ with $|t_0|$ sufficiently large,
$\mathcal{T}$ is a contraction on a suitable ball in
$C_\mu^{2,\alpha}((-\infty,t_0]\times S^{N-1})$.
Set
$$
\Gamma_{B,t_0}
=
\big\{
z \in C_\mu^{2,\alpha}((-\infty,t_0]\times S^{N-1})
:
\|z\|_{C_\mu^{2,\alpha}((-\infty,t_0]\times S^{N-1})} \le B
\big\}.
$$

We claim that for some fixed constant $B>0$ (independent of $\mu$) and for all
$t_0 < 0$ with $|t_0|$ sufficiently large, the  $\mathcal{T}$ maps
$\Gamma_{B,t_0}$ into itself. That means ${||\mathcal{T}(\phi)||}_{C_{\mu}^{2, \alpha}\left(\left(-\infty, t_0\right] \times S^{N-1}\right)} \leq B$ for
$||\phi||_{C_{\mu}^{2, \alpha}\left(\left(-\infty, t_0\right] \times S^{N-1}\right)} \leq B$.

First, by \eqref{64q} we have
$$
\|H(\widehat{z})\|_{C_\mu^1((-\infty,t_0]\times S^{N-1})}
\le C_1.
$$

Define
\begin{equation}
E(\phi)
=
(-p)
\int_0^1
\big[
(\widehat{z} + \Lambda + s\phi)^{-(p+1)}
-
\Lambda^{-(p+1)}
\big]\,
ds .
\tag{65}\label{qc}
\end{equation}
We get that $P(\phi) = \phi\, E(\phi)$.

Take any
$\phi \in C_\mu^{2,\alpha}((-\infty,t_0]\times S^{N-1}))$
satisfying
$\|\phi\|_{C_\mu^{2,\alpha}((-\infty,t_0]\times S^{N-1})} \le B$.
Note that
$$
|\widehat{z}| + |\nabla \widehat{z}| \le \varepsilon(t),
$$
where $\varepsilon(t)$ is a monotonically increasing function with
$\varepsilon(t)\to 0$ as $t\to -\infty$.
$$
|\phi| + |\nabla \phi|
\le
B e^{\mu t}.
$$

Hence, for all $t \le t_0$,
\begin{equation}
|E(\phi)| + |\nabla E(\phi)|
\le
C_2 \big( \varepsilon(t) + B e^{\mu t} \big),
\tag{66}\label{qv}
\end{equation}
therefore
$$
\begin{aligned}
\|P(\phi)\|_{C_\mu^1((-\infty,t_0]\times S^{N-1})}
\le
C_2 \big( \varepsilon(t_0) + B e^{\mu t_0} \big)
\|\phi\|_{C_\mu^1((-\infty,t_0]\times S^{N-1})}
\le
C_2 \big( \varepsilon(t_0) + B e^{\mu t_0} \big)\, B.
\end{aligned}
$$

By Theorem 3.7 we obtain
$$
\begin{aligned}
||\mathcal{T}(\phi)||_{C_\mu^{2, \alpha}\left(\left(-\infty, t_0\right] \times S^{N-1}\right)} & \leq C||-H(\hat{z})+P(\phi)||_{C_{\mu}^{0, \alpha}\left(\left(-\infty, t_0\right] \times S^{N-1}\right)} \\
& \leq C\left[C_1+C_2\left(\epsilon\left(t_0\right)+B e^{\mu t_0}\right) B\right],
\end{aligned}
$$
$C, C_1, C_2$ are all positive and independent of $t_0$.
Choose $B \ge 2 C C_1$ and  $|t_0|$ sufficiently large so that
$$
C C_2 \big( \varepsilon(t_0) + B e^{\mu t_0} \big)
\le \frac{1}{2}.
$$
It follows that
$$
\|\mathcal{T}(\phi)\|_{C_\mu^{2,\alpha}((-\infty,t_0] \times S^{N-1})}
\le B.
$$

This gives the required self-mapping property.
We now show that
$\mathcal{T} : \Gamma_{B,t_0} \to \Gamma_{B,t_0}$
is a contraction map; namely, for any
$\phi_1, \phi_2 \in \Gamma_{B,t_0}$, there exists
$\kappa \in (0,1)$ such that
\begin{equation}
\|\mathcal{T}(\phi_1) - \mathcal{T}(\phi_2)\|_{C_\mu^{2,\alpha}((-\infty,t_0] \times S^{N-1})}
\le
\kappa
\|\phi_1 - \phi_2\|_{C_\mu^{2,\alpha}((-\infty,t_0] \times S^{N-1})}.
\tag{67}\label{mm}
\end{equation}

Observe that
$$
\mathcal{T}(\phi_1) - \mathcal{T}(\phi_2)
=
\mathcal{L}^{-1}\big( P(\phi_1) - P(\phi_2) \big),
$$
and
$$
\begin{aligned}
P(\phi_1) - P(\phi_2)
&=
\phi_1 E(\phi_1) - \phi_2 E(\phi_2)
\\
&=
(\phi_1 - \phi_2) E(\phi_1)
\,+\,
\phi_2 \big( E(\phi_1) - E(\phi_2) \big).
\end{aligned}
$$

By \eqref{qc},
$$
E(\phi_1) - E(\phi_2)
=
(-p) \int_0^1
\big[
(\widehat{z} + \Lambda + s\phi_1)^{-(p+1)}
-
(\widehat{z} + \Lambda + s\phi_2)^{-(p+1)}
\big]\, ds.
$$

Thus,
$$
\left|E\left(\phi_1\right)-E\left(\phi_2\right)\right|+\left|\nabla\left(E\left(\phi_1\right)-E\left(\phi_2\right)\right)\right| \leq C\left(\left|\phi_1-\phi_2\right|+\left|\nabla\left(\phi_1-\phi_2\right)\right|\right) .
$$

For any $t \le t_0$ we have from \eqref{qv} that
$$
\begin{aligned}
&\big| P(\phi_1) - P(\phi_2) \big|
+
\big| \nabla ( P(\phi_1) - P(\phi_2) ) \big|
\\
&\qquad
\le
C \big( \varepsilon(t) + B e^{\mu t} \big)
\Big(
|\phi_1 - \phi_2|
+
|\nabla(\phi_1 - \phi_2)|
\Big),
\end{aligned}
$$
therefore
$$
\| P(\phi_1) - P(\phi_2) \|_{C_\mu^1((-\infty,t_0]\times S^{N-1})}
\le
C \big( \varepsilon(t_0) + B e^{\mu t_0} \big)
\| \phi_1 - \phi_2 \|_{C_\mu^1((-\infty,t_0]\times S^{N-1})}.
$$
By Theorem 3.7,
$$
\begin{aligned}
\| \mathcal{T}(\phi_1) - \mathcal{T}(\phi_2) \|_{C_\mu^{2,\alpha}((-\infty,t_0]\times S^{N-1})}
&\le
C \| P(\phi_1) - P(\phi_2) \|_{C_\mu^{0,\alpha}((-\infty,t_0]\times S^{N-1})}
\\
&\le
C \big( \varepsilon(t_0) + B e^{\mu t_0} \big)
\| \phi_1 - \phi_2 \|_{C_\mu^{2,\alpha}((-\infty,t_0]\times S^{N-1})}.
\end{aligned}
$$

Then we yield \eqref{mm} by taking $|t_0|$ sufficiently large.
By the contraction mapping principle, there exists
$\phi \in C_\mu^{2,\alpha}((-\infty,t_0]\times S^{N-1})$
such that
$\mathcal{T}(\phi) = \phi$.
This gives a solution $\phi$ to \eqref{64q}.
Consequently,
$z = \widehat{z} + \phi$
is a solution of $\mathcal{N}(z) = 0$ and satisfies \eqref{hh}.
\end{proof}
\begin{thm}
Let $\eta_{\star}(t) = d_F e^{\frac{2p - \alpha}{p+1} t}$ such that 
$\mathcal{L}_0(\eta_\star) = F e^{\frac{2p - \alpha}{p+1} t}$.

(1)Radial case.
Assume $F \neq 0$ and 
$\mu > \frac{(2p - \alpha) j}{p+1}$ 
for some integer $j \ge 2$, and that 
$\mu \notin \left\{ \frac{(2p - \alpha) k}{p+1} : k \ge 2 \right\}$.
Then there exists $t_0 < 0$ with $|t_0|$ sufficiently large, such that
there exists a smooth function $\widehat{\varphi}$ on $(-\infty,t_0]$ satisfying
$$
\widehat{z} = \eta_\star + \widehat{\varphi}(t),
$$
and $\widehat{z}$ satisfies the hypotheses of Theorem 4.8.

(2)Nonradial case.
Assume that the corresponding index sets $\mathcal{I}_\rho$ 
and $\mathcal{I}_{\widetilde{\rho}}$  given in \eqref{36.1} 
and \eqref{36.2} (with multiple number ignored) in the proof of 
Theorem 2.2 respectively.
$
\mu \notin \mathcal{I}_\rho \cup \mathcal{I}_{\widetilde{\rho}}
~\text{and}~
\mu > \rho_{j+4} > \frac{2p - \alpha}{p+1}
$
 in cases
\eqref{34}\eqref{35};
$
\mu > \rho_1
$in cases \eqref{30}\eqref{36}.
Let $\varphi$ be a solution of the equation 
$\mathcal{L}(\varphi) = 0$ on $\mathbb{R} \times \mathbb{S}^{N-1}$, and suppose that $\varphi(t,\theta) \to 0 ~~\text{as}~ t \to -
\infty ~\text{uniformly} ~\text{for}~ \theta \in{S^{N-1}}$.
Then there exists $t_0 < 0$ with $|t_0|$ sufficiently large, such that there
exists a smooth function $\widetilde{\varphi}$ on
$(-\infty,t_0]\times \mathbb{S}^{N-1}$ satisfying
$$
\widehat{z}
=
\eta_\ast + \varphi + \widetilde{\varphi},
$$
and $\widehat{z}$ satisfied the hypotheses of Theorem 4.7.
\end{thm}

\begin{proof}
We prove case (2); the proof of case (1) is similar and easier.

Take a function $\phi(t,\theta)$ such that
$\phi(t,\theta)\to 0$ uniformly for $\theta$ on $\mathbb{S}^{N-1}$
as $t\to -\infty$.
A direct computation yields
\begin{equation}
H(\eta_\ast + \phi)
=
\mathcal{L}(\phi)
-
\Big[
(\eta_\ast + \phi + \Lambda)^{-p}
-
\Lambda^{-p}
+
p \Lambda^{-(p+1)}(\eta_\ast + \phi)
\Big]
=
\mathcal{L}(\phi)
-
\sum_{i=2}^{\infty} b_i (\eta_\ast + \phi)^i.
\tag{68}\label{bf}
\end{equation}
(For convenience we write full infinite-series form of the nonlinear term ).

We now write the argument for case \eqref{34} with $k=1$ in the proof of Theorem 2.2.
The associated index sets $\mathcal{I}_\rho$ and $\mathcal{I}_{\widetilde{\rho}}$
are those given in \eqref{36.1} and \eqref{36.2} (ignoring multiple numbers).
In this case,
$$
\rho_1 = \sigma_1^{(1)},\quad
\rho_2 = \sigma_1^{(2)},\ \ldots,\
\rho_j = \sigma_1^{(j)},\quad
\rho_{j+1} = \frac{2p}{p+1},\quad
\rho_{j+2} = \sigma_1^{(j+1)},\ldots.
$$
Let $K \ge j+4$ be the largest integer such that $\rho_K < \mu$,
and let $\widetilde{K}$ be the largest integer such that
$\widetilde{\rho}_{\widetilde{K}} < \mu$.
Since the kernel of $\mathcal{L}_0$ contains no function
to zero as $t\to -\infty$. And  the kernel of $\mathcal{L}_k$
contains one exponentially decaying function
$\psi_k^{+}(t)=e^{\sigma_1^{(k)} t}$ and one exponentially growing function
$\psi_k^{-}(t)=e^{\sigma_2^{(k)} t}$ as  . Without loss of generality,we
assume that the solution $\varphi$ of $\mathcal{L}(\varphi)=0$ has the form
\begin{equation}
\varphi(t,\theta)
=
\sum_{i=1}^j c_i Q_i(\theta) e^{\rho_i t}
+
\sum_{i=j+2}^K c_i Q_{\,i-1}(\theta) e^{\rho_i t},
\tag{69}\label{gf}
\end{equation}
where $c_i$ are constants.
This is because any term of the form $e^{\rho_i t}$ with $i > K$ appearing in
$\varphi$ would contribute only terms $e^{\widetilde{\rho}_\ell t}$ with
$\widetilde{\rho}_\ell > \mu$ in $H(\eta_\ast + \varphi)$.

We first consider the case that
\begin{equation}
\mathcal{I}_\rho \cap \mathcal{I}_{\tilde{\rho}}=\emptyset
\tag{70}\label{JK}
\end{equation}
We shall prove that it can successively construct
$\widetilde{\varphi}_0, \widetilde{\varphi}_1, \ldots, \widetilde{\varphi}_K$
such that for every $l = 0,1,\ldots,\widetilde{K}$,
\begin{equation}
H\big(\eta_\ast + \varphi
+ \widetilde{\varphi}_0 + \cdots + \widetilde{\varphi}_l\big)
=
O\big( e^{\widetilde{\rho}_{l+1} t} \big).
\tag{71}\label{gg}
\end{equation}

We first take $\phi = \varphi$, with $\varphi$ given by \eqref{gf}.
Then by \eqref{bf} and the fact that $\mathcal{L}(\varphi) = 0$, we have
$$
H(\eta_\ast + \varphi)
=
\sum_{n_1 + \cdots + n_{r_1} \ge 2}
a_{n_1 \ldots n_{r_1}}
e^{(n_1 \rho_1 + \cdots + n_{r_1} \rho_{r_1}) t}
Q_1^{n_1} \cdots Q_j^{n_j} Q_0^{n_{j+1}}
Q_{j+1}^{n_{j+2}} \cdots Q_{r_1-1}^{n_{r_1}},
$$
where $n_1,\ldots,n_{r_1}$ are nonnegative integers and
$a_{n_1 \ldots n_{r_1}}$ are constants.
By the definition of $\mathcal{I}_{\widetilde{\rho}}$ that each
$n_1 \rho_1 + \cdots + n_{r_1} \rho_{r_1}$ is equal to some $\widetilde{\rho}_i$.
Hence
\begin{equation}
H(\eta_\ast + \varphi)
=
\sum_{i=1}^{\widetilde{K}}
\left\{
\sum_{m=0}^{\widetilde{M}_i}
a_{im} Q_m(\theta)
\right\}
e^{\widetilde{\rho}_i t}
+ O\big( e^{\widetilde{\rho}_{K+1} t} \big),
\tag{72}\label{ss}
\end{equation}
where $\widetilde{M}_i$ is given by \eqref{30q} and $a_{im}$ are constants.
In particular,
$$
H\left(\eta_*+\varphi\right)=O\left(e^{\tilde{\rho}_1 t}\right),
$$
Here $\widetilde{\rho}_1 = 2\rho_1$. Hence  \eqref{gg} holds when $l = 0$ and
$\widetilde{\varphi}_0 = 0$.

Assume that we have already constructed
$\widetilde{\varphi}_0, \widetilde{\varphi}_1, \ldots, \widetilde{\varphi}_{l-1}$
such that \eqref{gg} holds for $0,1,\ldots,l-1$.  
We now consider the case $l$.  
Set
\begin{equation}
\widetilde{\varphi}_l(t,\theta)
=
\left( \sum_{m=0}^{\widetilde{M}l} c_{lm} Q_m(\theta) \right)
e^{\widetilde{\rho}_l t},
\tag{73}\label{ll}
\end{equation}
where $c_{lm}$ are constants to be determined.
A computation similar to that leading to \eqref{ss} yields
$$
\begin{aligned}
H\big(\eta_\ast + \varphi
+ \widetilde{\varphi}_0 + \cdots + \widetilde{\varphi}_l\big)
=
\mathcal{L}(\widetilde{\varphi}_1)
+ \cdots
+ \mathcal{L}(\widetilde{\varphi}_l)
+
\sum_{i=1}^{\widetilde{K}}
\left\{
\sum_{m=0}^{\widetilde{M}_i}
a_{im} Q_m(\theta)
\right\}
e^{\widetilde{\rho}_i t}
+ O\big( e^{\widetilde{\rho}_{\widetilde{K}+1} t} \big),
\end{aligned}
$$
where $a_{im}$ are constants whose values may differ from those in (158).
By the induction hypothesis,
$$
H\big(\eta_\ast + \varphi
+ \widetilde{\varphi}_0 + \cdots + \widetilde{\varphi}_{l-1}\big)
=
\sum_{i=l}^{\widetilde{K}}
\left\{
\sum_{m=0}^{\widetilde{M}_i}
a_{im} Q_m(\theta)
\right\}
e^{\widetilde{\rho}_i t}
+ O\big( e^{\widetilde{\rho}_{\widetilde{K}+1} t} \big),
$$
and hence
$$
H\big(\eta_\ast + \varphi
+ \widetilde{\varphi}_0 + \cdots + \widetilde{\varphi}_l\big)
=
\mathcal{L}(\widetilde{\varphi}_l)
+
\sum_{i=l}^{\widetilde{K}}
\left\{
\sum_{m=0}^{\widetilde{M}_i}
a_{im} Q_m(\theta)
\right\}
e^{\widetilde{\rho}_i t}
+ O\big( e^{\widetilde{\rho}_{\widetilde{K}+1} t} \big).
$$

Note that $\widetilde{\varphi}_l$ does not contribute to the coefficients
$a_{lm}$. We choose $\widetilde{\varphi}_l$ so that
$$
\mathcal{L}\left(\tilde{\varphi}_{l}\right)
=-\left\{\sum_{m=0}^{\tilde{M}_l} a_{l m} Q_m(\theta)\right\} e^{\tilde{\rho}_l t} .
$$

In this way we obtain \eqref{gg} for the index $l$.
With $\widetilde{\varphi}_l$ given by \eqref{ll}, it remains only to solve
\begin{equation}
\mathcal{L}_m\big( c_{lm} e^{\widetilde{\rho}_l t} \big)
= - a_{lm} e^{\widetilde{\rho}_l t}
\tag{74}\label{aa}
\end{equation}
for $m = 0,1,\ldots,\widetilde{M}_l$.

Since $\rho_m \ne \widetilde{\rho}_l$ for all $m,l$, we can find constants
$c_{lm}$ satisfying \eqref{aa}.
Indeed, we can write an explicit formula for $c_{lm} e^{\widetilde{\rho}_l t}$ in terms of
$a_{lm} e^{\widetilde{\rho}_l t}$ by using the basis of $\operatorname{Ker}(\mathcal{L}_m)$ given in Lemma 4.3.

if $0 < \rho_m < \widetilde{\rho}_l$, this representation is provided by \eqref{61};
if $\rho_m > \widetilde{\rho}_l$, we simply replace the first integral in \eqref{61}
with the one from $t$ to $t_0$.
This completes the induction.

In conclusion, we set
\begin{equation}
\widetilde{\varphi}(t,\theta)
=
\sum_{i=1}^{\widetilde{K}}
\left\{
\sum_{m=0}^{\widetilde{M}_i}
c_{im} Q_m(\theta)
\right\}
e^{\widetilde{\rho}_i t},
\tag{75}\label{1Q}
\end{equation}
where $c_{im}$ are constants. Then
$$
H(\eta_\ast + \varphi + \widetilde{\varphi})
=
O\big( e^{\widetilde{\rho}_{\widetilde{K}+1} t} \big)
=
O(e^{\mu t}).
$$

An similar estimate holds for the gradient of
$H(\eta_\ast + \varphi + \widetilde{\varphi})$.
This completes the proof in this case.

We now consider the general situation, in which some $\rho_i$ can be written
as a positive--integer linear combination of $\rho_1,\ldots,\rho_{i-1}$.
We briefly indicate how to modify the above argument to handle this case.
The modification mainly concerns \eqref{aa}.
If for some coefficient $a_{lm}$ we have $\rho_m = \widetilde{\rho}_l$,
then instead of choosing only a constant $c_{lm}$ as in \eqref{aa},
we can find constants $c_{l0m}$ and $c_{l1m}$ such that
$$
\mathcal{L}_m\left(\left(c_{l 0 m}+t c_{l 1 m}\right) e^{\tilde{\rho}_l t}\right)
=-a_{l m} e^{\tilde{\rho}_l t} .
$$

Such powers of $t$ will generate higher-order powers of $t$ in the process of the iteration.
In general, for a nonnegative integer $J$, if constants $a_{ljm}$ are given for
$j = 0,1,\ldots,J$, then we can find constants $c_{ljm}$, $j = 0,1,\ldots,J+1$, such that
$$
\mathcal{L}_m\left( \sum_{j=0}^{J+1} c_{ljm} t^j e^{\widetilde{\rho}_l t} \right)
=
\sum_{j=0}^{J} a_{ljm} t^j e^{\widetilde{\rho}_l t}.
$$

Therefore, in the general case we do not adopt \eqref{1Q}, but instead take
$$
\widetilde{\varphi}(t,\theta)
=
\sum_{i=1}^{\widetilde{K}}
\sum_{j=0}^{i}
\left\{
\sum_{m=0}^{\widetilde{M}_i}
c_{ijm} Q_m(\theta)
\right\}
t^j e^{\widetilde{\rho}_i t},
$$
where $c_{ijm}$ are constants.
This completes the proof of the proposition.
\end{proof}

In particular, we have now found a function $z(t,\theta)$ satisfying the required conditions.
Recalling that
$$
u(x) = u_s(|x|) + |x|^{\frac{\alpha+2}{p+1}} z(\ln |x|,\theta),
$$
we obtain the desired solution $u$, and the proof of Theorem 1.1 and Theorem 1.2 is thereby completed.

\section*{Acknowledgments}
 Y. Zhang is sponsored by NSFC [No. 12271505] and STCSM [No. 22DZ2229014].


\begin{thebibliography}{1}
\bibitem{11z}
Juan Davila, Dong Ye. On finite Morse index solutions of two equations with negative exponent[J]. Proceedings of the Royal Society of Edinburgh, 2013, 143A, 121-128.

\bibitem{1z}
Z. M. Guo, X. Huang, and F. Zhou, "Radial symmetry of entire solutions of a bi-harmonic equation with exponential nonlinearity," J. Funct. Anal. 268, 1972--2004
(2015).

\bibitem{3z}
H. Zou, "Symmetry of positive solutions of $\Delta u + u^p = 0$ in $R^n$
," J. Differ. Equations 120, 46--88 (1995).

\bibitem{4z}
Q. Han, X. X. Li, and Y. C. Li, "Asymptotic expansions of solutions of the Yamabe equation and the $\sigma_k$-Yamabe equation near isolated singular points," Commun. Pure
Appl. Math. 74, 1915--1970 (2021).

\bibitem{5z}
Stein, E. M.; Weiss, G. Introduction to Fourier analysis on Euclidean spaces. Princeton Uni- versity Press, Princeton, N.J., (1971)

\bibitem{6z}
Hongxia Guo, Zongming Guo, Ke Li. Positive solutions of a semilinear elliptic equation with singular nonlinearity[J]. Journal of Mathematical Analysis and Applications,  323-344,(2006)

\bibitem{2z}
Marie-Francoise Bidaut-Veron, Victor Galaktionov, Philippe Grillot, Laurent Veron. Singularities for a 2-Dimensional Semilinear Elliptic Equation with a Non-Lipschitz Nonlinearity[J]. Journal of Differential Equations, 1999, 154(2), 318-338.

\bibitem{7z}
Yujin Guo, Yanyan Zhang, Feng Zhou. Singular behavior of an electrostatic-elastic membrane system with an external pressure[J]. Nonlinear Analysis,  190. 111611.(2020)

\bibitem{8z}
Z. M. Guo; F. S. Wan. Positive rupture solutions of steady states for thin-film-type equations. J. Math. Phys. 65, 041510 (2024)

\bibitem{9z}
 Qing Li and Yanyan Zhang, Study on the behaviors of rupture solutions for a class of elliptic MEMS equations in  . Journal of Differential Equations 414  817--841 (2025)

\bibitem{21z} 
Juan Dávila, Kelei Wang, and Juncheng Wei, Qualitative analysis of rupture solutions for a MEMS problem, Ann. Inst. H. Poincaré C Anal. Non Linéaire 33(1) (2016), 221–242.

\bibitem{10z}
Xuyan Chen, Hiroshi Matano, Laurent Veron. Anisotropic singularities of solutions of nonlinear elliptic equations in   [J]. Journal of Functional Analysis, 1, 50-97,(1989)

\bibitem{12z}
Y. H. Du and Z. M. Guo, "Positive solutions of an elliptic equation with negative exponent: Stability and critical power," J. Differ. Equations 246, 2387--2414 (2009).

\bibitem{19z} 
Zongming Guo and Juncheng Wei, Rupture solutions of an elliptic equation with a singular nonlinearity, Proc. R. Soc. Edinb. Sect. A Math. 144(5) (2014), 905–924. 

\bibitem{13z}
Q. Han and Y. C. Li, "Singular solutions to the Yamabe equation with prescribed asymptotics," J. Differ. Equations 274, 127--150 (2021).

\bibitem{14z}
Han, Z.-C.; Li, Y.; Teixeira, E. V. Asymptotic behavior of solutions to the k-Yamabe equation
near isolated singularities. Invent. Math. 182 (2010), no. 3, 635--684.

\bibitem{15z}
Korevaar, N.; Mazzeo, R.; Pacard F.; Schoen, R. Refined asymptotics for constant scalar
curvature metrics with isolated singularities. Invent. Math. 135 (1999).

\bibitem{22z} 
H. Chen, Y. Wang and F. Zhou, On semilinear elliptic equation with negative exponent arising from a closed MEMS model, Z. Angew. Math. Phys. 75 (2024), Article 3, 1–25

\bibitem{16z}
Z. M. Guo, J. Y. Li, and F. S. Wan, "Asymptotic behavior at the isolated singularities of solutions of some equations on singular manifolds with conical metrics,"
Commun. Partial Differ. Equations 45, 1647--1681 (2020).

\bibitem{17z}
Mazzeo, R.; Pacard, F. Constant scalar curvature metrics with isolated singularities. Duke Math.
J. 99 (1999), no. 3, 353--418. 

\bibitem{18z} 
Mazzeo, R.; Pollack, D.; Uhlenbeck, K. Moduli spaces of singular Yamabe metrics. J. Amer.
Math. Soc. 9 (1996), no. 2, 303--344.

\bibitem{20z} 
Wei Wang and Zhifei Zhang, Fine structure of rupture set for semilinear elliptic equation with singular nonlinearity, arXiv:2411.16048v2 [math.AP] (2024). 
 
\bibitem{23z}
J.R. Beckham, J.A. Pelesko, An electrostatic-elastic membrane system with an external pressure,Math. Comput.
Modelling 54 (2011) 2551–3212.

\bibitem{24z}
P. Esposito, N. Ghoussoub, Y. Guo, Mathematical analysis of partial differential equations modeling electrostatic MEMS,
Courant Lect. Notes Math. 20 (2010) 318.

\bibitem{25z}
P.Laurencot, Philippe, C.Walker, Some singular equations modeling MEMS, Bull. Amer. Math. Soc. 54 (2017) 437–479.
\end{thebibliography}
\end{document}